\documentclass{amsart}
\usepackage{amsthm,amsmath,amssymb,amscd,graphics,psfig,epic,eepic}

\newcommand{\entrylabel}[1]{\mbox{\textsf{{\rm c}1}}\hfil}
{\end{list}}
{
   \newtheorem{theorem}{Theorem}[subsection]                     
   \newtheorem{proposition}[theorem]{Proposition}     
   \newtheorem{lemma}[theorem]{Lemma}

}
{\theoremstyle{definition}
   
   \newtheorem{example}[theorem]{Example}
   \newtheorem{definition}[theorem]{Definition}
}
{\theoremstyle{remark}
   \newtheorem*{remark}{Remark}
   \newtheorem*{remarks}{Remarks}
}
\newcommand{\Hi}{\cite{Hir}}

\newcommand{\Hii}{\cite{Hir1}}
\newcommand{\Vi}{\cite{Vi}}

\newcommand{\BM}{\cite{BM}}

\newcommand{\cC}{{\mathcal C}}
\newcommand{\cD}{{\mathcal D}}

\newcommand{\cH}{{\mathcal H}}
\newcommand{\cI}{{\mathcal I}}
\newcommand{\cJ}{{\mathcal J}}

\newcommand{\cM}{{\mathcal M}}
\newcommand{\cN}{{\mathcal N}}
\newcommand{\cO}{{\mathcal O}}

\newcommand{\supp}{{\operatorname{supp}}}

\newcommand{\res}{{\operatorname{res}}}

\newcommand{\Sub}{\operatorname{Sub}}

\newcommand{\Ga}{\operatorname{Gal}}

\newcommand{\Spec}{\operatorname{Spec}}

\newcommand{\Der}{{\operatorname{Der}}}

\newcommand{\ord}{{\operatorname{ord}}}

\newcommand{\codim}{\operatorname{codim}}

\newcommand{\inv}{{\operatorname{inv}}}

\newcommand{\Reg}{{\operatorname{Reg}}}

\begin{document}


\title[Simple Hironaka resolution in  characteristic zero]{Simple Hironaka resolution in  characteristic zero}

\author{Jaros{\l}aw W{\l}odarczyk}
\thanks{The author was supported in part by NSF grant DMS-0100598 and Polish KBN
grant 2 P03 A 005 16}
\address{Department of Mathematics\\Purdue University\\West
Lafayette, IN-47907\\USA}

\email{wlodar@math.purdue.edu, jwlodar@mimuw.edu.pl}
\date{\today}
\begin{abstract}
Building upon work  of 
Villamayor and Bierstone-Milman we give a
 proof of the canonical Hironaka principalization and desingularization.
 The idea of "homogenized ideals" introduced in the paper gives {\it a priori} the canonicity of algorithm and radically simplifies the proof.

\end{abstract}
\maketitle

\tableofcontents
\addtocounter{section}{-1}

\section{Introduction} In the present paper we give a short proof of the Hironaka theorem on resolution of singularities.  Recall that in the classical approach to the problem of embedded resolution originated by Hironaka (see \cite{Hir}) and later developed and simplified  by Bierstone-Milman  (see \cite{BM}) and Villamayor (see \cite{Vi}) an invariant which plays the role of a measure of singularities is constructed. 
The invariant is upper semicontinuous and defines a stratification of the ambient space.
This invariant drops after the blow-up of the maximal stratum. 
It determines the centers of the resolution and allow to  patch up local desingularizations to a global one. Such an invariant carries 
a rich information about singularities and the resolution process. 
The definition of the invariant is quite involved. What adds to complexity is that the invariant is defined within some rich inductive scheme encoding the  desingularization and assuring its canonicity ( Bierstone-Milman's towers of local blow-ups with {\it admissible centers} and Villamayor's {\it generalized basic objects}) 
(see also Encinas-Hauser \cite{EH}).

The idea of forming the invariant is based upon the observation due to H. Hironaka that the resolution process
controlled by the order or Hilbert-Samuel function can be reduced to the resolution process on some smooth hypersurface, called a {\it hypersurface of maximal contact} (see \cite{Hir},).
The reduction to a hypersurface of maximal contact is not canonical and for two different hypersurfaces of maximal contact we get two different objects loosely related but having the same invariant. 
To make all the process canonical and relate the objects induced by restrictions we either interpret the invariant in a canonical though quite technical way (so called "Hironaka trick")  or build a relevant canonical resolution datum  (\cite{BM1},\cite{BM2},\cite{BM3},\cite{EH},\cite{Ha},\cite{Hir},\cite{Vi},\cite{Vi1},\cite{Vi2}). (Both ways are strictly related). 

The approach we propose in this paper is
based upon the above mentioned reduction procedure and two simple observations.
\begin{enumerate}
\item The resolution process defined as a sequence of the suitable blow-ups of ambient spaces can be applied simultaneously not only to the given singularities but rather to a class of equivalent singularities obtained by simple arithmetical modifications. This means that we can "tune" singularities before resolving them.  
\item In the equivalence class we can choose a convenient representative given by the {\it homogenized  ideals} introduced in the paper. The restrictions of  homogenized ideals to different hypersurfaces of maximal contact define locally analytically isomorphic singularities. Moreover  the  local isomorphism of 
hypersurfaces of  maximal contact is defined by a local analytic automorphism of the ambient space preserving all
the relevant resolutions.
\end{enumerate}

As a consequence neither constructing painstaking structures nor studying elaborate invariants are required for the proof of the
resolution theorems. All we need is the existence of a canonical functorial resolution in lower dimensions. The use of  an invariant controlling the whole resolution procedure  is now not essential and can be completely avoided. Despite that, we  introduce the invariant in the process of constructing the algorithm mainly to mark the progress towards the desingularization.  

The strategy of the proof we formulate here is essentially the same as the one found by Hironaka and simplified by Bierstone-Milman and Villamayor.

Although the preliminary set-up (with some modifications and refinements) and general strategy  of the algorithm build upon work of Villamayor (see \cite{Vi},  \cite{Vi1}, \cite{Vi2}) the present proof leads to essentially the same invariant  as in  Bierstone-Milman papers (see \cite{BM}, \cite{BM1}, \cite{BM2})

The presented proof is fairly elementary, constructive and self-contained.

The methods used in this paper can be applied to desingularization of analytic space;   we deal with the analytic case in a separate paper. 

\section{Formulation of the main theorems}
All algebraic varieties in this paper are defined over a  ground field of characteristic zero. The assumption of characteristic zero is only needed for the local existence of a hypersurface of  maximal contact (Lemma \ref{le: Gi}). 

We give a proof of the following Hironaka Theorems (see \cite{Hir}):
\begin{enumerate}

\item {\bf Canonical Principalization}
\begin{theorem} \label{th: 1} Let ${\cI}$ be a sheaf of ideals on a smooth algebraic variety $X$.
There exists a principalization of ${\cI}$, that is, a sequence 

$$ X=X_0 \buildrel \sigma_1 \over\longleftarrow X_1 
\buildrel \sigma_2 \over\longleftarrow X_2\longleftarrow\ldots
\longleftarrow X_i \longleftarrow\ldots \longleftarrow X_r =\widetilde{X}$$

of blow-ups $\sigma_i:X_{i-1}\leftarrow X_{i}$ of smooth centers $C_{i-1}\subset
 X_{i-1}$
such that

\begin{enumerate}
\item The exceptional divisor $E_i$ of the induced morphism $\sigma^i=\sigma_1\circ \ldots\circ\sigma_i:X_i\to X$ has only  simple normal 
crossings and $C_i$ has simple normal crossings with $E_i$.

\item The
total transform $\sigma^{r*}({\cI})$ is the ideal of a simple normal
crossing divisor
$\widetilde{E}$ which is  a natural  combination of the irreducible components of the divisor ${E_r}$.

\end{enumerate}
The morphism $(\widetilde{X},\widetilde{\cI})\rightarrow(X,{\cI}) $ defined by the above principalization  commutes with smooth morphisms and  embeddings of ambient varieties. It is equivariant with respect to any group action not necessarily preserving the ground field $K$.

\end{theorem}

\item {\bf  Weak-Strong Hironaka Embedded  Desingularization}
\begin{theorem} \label{th: emde} \label{th: 2} Let $Y$ be a subvariety of a smooth variety
$X$ over a field of characteristic zero.
There exists a  sequence $$ X_0=X \buildrel \sigma_1 \over\longleftarrow X_1 
\buildrel \sigma_2 \over\longleftarrow X_2\longleftarrow\ldots
\longleftarrow X_i \longleftarrow\ldots \longleftarrow X_r=\widetilde{X}$$ of
blow-ups  $\sigma_i:X_{i-1}\longleftarrow X_{i}$ of smooth centers $C_{i-1}\subset
 X_{i-1}$ such that

\begin{enumerate}

\item The exceptional divisor $E_i$ of the induced morphism $\sigma^i=\sigma_1\circ \ldots\circ\sigma_i:X_i\to X$ has only  simple normal 
crossings and $C_i$ has simple normal crossings with $E_i$.

\item Let $Y_i\subset X_i$ be the strict transform of $Y$. All centers $C_i$ are disjoint from the set $\Reg(Y)\subset Y_i$ of points where   $Y$  (not $Y_i$) is smooth
(and are not necessarily contained in $Y_i$).

\item  The strict transform $\widetilde{Y}:=Y_r$ of  $Y$  is smooth and 
has only simple normal crossings with the exceptional divisor $E_r$.

\item The morphism $(X,{Y})\leftarrow (\widetilde{X},\widetilde{Y})$ defined by the embedded desingularization commutes with smooth morphisms and  embeddings of ambient varieties. It is equivariant with respect to any group action not necessarily preserving $K$.

\item (Strengthening of Bravo-Villamayor) (see \cite{BV}) $$\sigma^*(\cI_Y)=\cI_{\widetilde{Y}}\cI_{\widetilde{E}},$$ where $\cI_{\widetilde{Y}}$ is the sheaf of ideals of the subvariety $\widetilde{Y}\subset\widetilde{X}$ and $\cI_{\widetilde{E}}$ is the sheaf of ideals of a simple normal
crossing divisor
$\widetilde{E}$ which is  a natural  combination of the irreducible components of the divisor ${E_r}$.

 \end{enumerate}

\end{theorem}

\item {\bf Canonical Resolution of Singularities}
\begin{theorem} \label{th: 3} Let $Y$ be an algebraic variety over a field of characteristic zero.

There exists a canonical desingularization of $Y$ that is
a smooth variety $\widetilde{Y}$ together with a proper birational morphism $\res_Y: \widetilde{Y}\to Y$ which is functorial with respect to smooth morphisms. 
For any smooth morphism $\phi: Y'\to Y$ there is a natural lifting $\widetilde{\phi}: \widetilde{Y'}\to\widetilde{Y}$ which is a smooth morphism.

In particular $\res_Y: \widetilde{Y}\to Y$ is  an isomorphism over the nonsingular part of $Y$. Moreover $\res_Y$ is equivariant with respect to any group action not necessarily preserving the ground field.
\end{theorem}
\end{enumerate}

\section{Preliminaries}
To simplify our considerations we shall assume that the ground field is algebraically closed. 
At the end of the paper we deduce the theorem for an arbitrary ground field of characteristic zero.

\subsection{Resolution of marked ideals} 

For any sheaf of ideals ${\cI}$ on  a smooth  variety $X$ and any point $x\in X$ we denote by
$$\ord_x({\cI}):=\max\{i\mid \cI\subset m_x^i\}$$ the {\it order} of ${\cI}$ at $x$. 
(Here $m_x$ denotes the maximal ideal of $x$.)

\begin{definition}(Hironaka (see \Hi, \Hii),  Bierstone-Milman (see \BM),Villamayor (see \Vi)) 
A {\it marked ideal} (originally a {\it basic object}  of  Villamayor) is a collection $(X,{\cI},E,\mu)$, where $X$ is a smooth variety,
${\cI}$ is a  sheaf of ideals on $X$, $\mu$ is a nonnegative integer and $E$ is a totally ordered collection of  divisors 
whose  irreducible components are pairwise disjoint and all have multiplicity one. Moreover the irreducible components of divisors in $E$ have simultaneously 
simple normal crossings. \end{definition}
\begin{definition}(Hironaka (\Hi, \Hii), Bierstone-Milman (see \BM),Villamayor (see \Vi))
By the {\it support} (originally {\it singular locus}) of $(X,{\cI},E,\mu)$ we mean
$$\supp(X,{\cI},E,\mu):=\{x\in X\mid \ord_x(\cI)\geq \mu\}.$$
\end{definition}

\begin{remarks} 
\begin{enumerate}
\item The ideals with assigned orders or functions with assigned multiplicities and their supports are key objects in proofs of Hironaka, Villamayor and Bierstone-Milman (see \cite{Hir}. 
Hironaka introduced the notion of {\it idealistic exponent}. Then various modifications 
of this definition were considered in the papers of Bierstone-Milman ({\it presentation of invariant}) and Villamayor ( {\it basic objects}).
In our proof we stick to   Villamayor's presentation of his basic objects (and their resolutions). Our marked ideals
are essentially the same notion as basic objects. However because of some technical differences and in order to introduce more suggestive terminology we shall call them  marked ideals.
\item Sometimes for simplicity we shall represent marked ideals $(X,{\cI},E,\mu)$ as couples $({\cI},\mu)$ or even ideals ${\cI}$.

\item For any sheaf of ideals ${\cI}$ on $X$ we have
 $\supp({\cI},1)=\supp({\cI})$.
\item For any marked ideals $({\cI},\mu)$ on $X$,
 $\supp({\cI},\mu)$ is a closed subset of $X$ (Lemma \ref{le: Vi1}).
\end{enumerate}
\end{remarks}

\begin{definition}(Hironaka (see \Hi, \Hii),  Bierstone-Milman (see \BM),Villamayor (see \Vi)) 
By a {\it resolution} of $(X,{\cI},E,\mu)$ we mean
a sequence of blow-ups $\sigma_i:X_i\to X_{i-1}$ of disjoint unions of smooth centers $C_{i-1}\subset
 X_{i-1}$, 

$$ 
X_0=X\buildrel \sigma_1 \over\longleftarrow X_1 
\buildrel \sigma_2 \over\longleftarrow X_2 \buildrel
\sigma_3 \over \longleftarrow\ldots
X_i\longleftarrow \ldots \buildrel \sigma_{r}  \over\longleftarrow X_r,$$

\noindent which defines a sequence of  marked ideals
$(X_i,{\cI}_i,E_i,\mu)$ where 

\begin{enumerate}
\item $C_i
\subset
\supp(X_i,{\cI}_i,E_i,\mu)$. 
\item $C_i$ has simple normal crossings with $E_i$.
\item ${\cI}_i=\cI(D_i)^{-\mu}\sigma_i^*({\cI}_{i-1})$, where $\cI(D_i)$ is the ideal of  the  exceptional divisor $D_i$
of $\sigma_i$.


\item $E_i=\sigma_i^{\rm c}(E_{i-1})\cup \{D_i\}$, where  $\sigma_i^{\rm c}(E_{i-1})$ is the set of strict transforms of divisors in $E_{i-1}$.
\item The order  on $\sigma_i^{\rm c}(E_{i-1})$ is defined by the order on $E_{i-1}$ while $D_i$ is the maximal element of $E_i$.
\item $\supp(X_r,{\cI}_r,E_r,\mu)=\emptyset$.
\end{enumerate}

\begin{definition}
The   sequence of morphisms which are either isomorphisms or  blow-ups satisfying  conditions (1)-(5) is called a {\it multiple test  blow-up}. The number of morphisms in a multiple test blow-up will be called its {\it length}.
\end{definition}
\begin{definition}
An {\it extension} of a multiple test blow-up (or a resolution) $(X_i)_{0\leq i\leq m}$ is a sequence $(X'_j)_{0\leq j\leq m'}$ of blow-ups and isomorphisms
$X'_{0}=X'_{j_0} =\ldots=X'_{j_1-1}\leftarrow X'_{j_1}=\ldots=X'_{j_2-1}\leftarrow \ldots X'_{j_m}=\ldots=X'_{m'}$,
where $X'_{j_i}=X_i$.
\end{definition}

\begin{remarks} 
\begin{enumerate}
\item The definition of  extension  arises naturally when we pass to open subsets of the considered ambient variety $X$.
\item The notion of a {\it multiple test  blow-up} is analogous to the notions of {\it test } or {\it admissible} blow-ups considered by Hironaka, Bierstone-Milman and Villamayor.

\end{enumerate}
\end{remarks}

\subsection{Transforms of marked ideal and controlled transforms of functions} 
In the setting of the above definition we shall call $$({\cI}_i,\mu):=\sigma_i^{{\rm c}}({\cI}_{i-1},\mu)$$ 
\end{definition}
\noindent a {\it transform of the marked ideal} or  {\it controlled transform} of $(\cI,\mu)$. It makes sense for a single blow-up in a multiple test  blow-up as well as for a multiple test  blow-up. Let  $\sigma^i:=\sigma_1\circ\ldots\circ\sigma_i: X_i\to X$ be a composition of consecutive morphisms of a multiple test  blow-up. Then in the above setting
$$({\cI}_i,\mu)=\sigma^{i{\rm c}}({\cI},\mu).$$ 
We shall also denote the controlled transform $\sigma^{i{\rm c}}({\cI},\mu)$ by
$(\cI,\mu)_i$ or $[\cI,\mu]_i.$

The controlled transform can also be defined for local sections $f\in\cI(U)$. Let $\sigma: X\leftarrow X'$ be a blow-up of a smooth center $C\subset \supp(\cI,\mu)$ defining transformation of marked ideals $\sigma^{\rm c}(\cI,\mu)=(\cI',\mu)$. Let $f\in \cI(U)$ be a section of a sheaf of ideals. Let $U'\subseteq\sigma^{-1}(U)$ be an open subset for which the sheaf of ideals of the exceptional divisor is generated by a function $y$. The function $$g=y^{-\mu}(f\circ\sigma) \in  \cI(U')$$ \noindent 
is a {\it controlled transform} of $f$ on $U'$ (defined up to an invertible function). As before we extend it  to any multiple test  blow-up.

The following lemma shows that the notion of controlled transform is  well defined.
\begin{lemma} Let $C\subset \supp({\cI},\mu)$ be a smooth center of the blow-up $\sigma: X\leftarrow X'$ and let $D$ denote the exceptional divisor.  Let ${\cI}_C$ denote the sheaf of ideals defined by $C$. Then
\begin{enumerate}
\item $\cI \subset {\cI}_C^\mu$.
\item  $\sigma^*({\cI})\subset ({\cI}_D)^\mu$.
\end{enumerate}
\end{lemma}
\noindent {\bf Proof.} (1) We can assume that the ambient variety $X$ is affine. Let $u_1,\ldots,u_k$ be parameters generating ${\cI}_C$
Suppose $f \in \cI \setminus {\cI}_C^\mu$. Then we can write $f=\sum_{\alpha}c_\alpha u^\alpha$, where either $|\alpha|\geq \mu$ or  $|\alpha|< \mu$ and $c_\alpha\not\in {\cI}_C$.
By the assumption there is $\alpha$ with $|\alpha|< \mu$ such that $c_\alpha\not \in {\cI}_C$. Take  $\alpha$  with  the smallest $|\alpha|$. There is a point $x\in C$ for which $c_\alpha(x)\neq 0$ and in the Taylor expansion of $f$ at $x$ there is a term $c_\alpha(x) u^\alpha$. Thus $\ord_x(\cI)< \mu$. This contradicts to the assumption $C\subset \supp({\cI},\mu)$. 

(2) $\sigma^*({\cI})\subset \sigma^*({\cI}_C)^\mu=({\cI}_D)^\mu$.\qed
\subsection{Hironaka resolution principle}

Our proof is based upon the following principle which can be traced back to Hironaka and was used by Villamayor in his simplification of Hironaka's algorithm:
\begin{eqnarray}  &\textbf{  (Canonical) Resolution of marked ideals $(X,\cI,E,\mu)$}\\ &\textbf{$\Downarrow$}\nonumber\\
& \textbf{ (Canonical) Principalization of the sheaves $\cI$ on $X$}\\ 
&\textbf{$\Downarrow$}\nonumber\\
&\textbf{(Canonical) Weak Embedded
Desingularization of subvarieties $Y\subset X$}\\  
&\textbf{$\Downarrow$}\nonumber\\
&\textbf{(Canonical) Desingularization} \end{eqnarray}

(1)$\Rightarrow$(2)  It follows immediately from the definition that a resolution of $(X,{\cI},\emptyset,1)$
determines a principalization of ${\cI}$. Denote by $\sigma: X\leftarrow \widetilde{X}$ the morphism defined by a resolution of $(X,\cI,\emptyset,1)$. The controlled transform $(\widetilde{\cI},1):=\sigma^{\rm c}(\cI,1)$ has the empty support. Consequently, $V(\widetilde{\cI})=\emptyset$,  and thus $\widetilde{\cI}$ is equal to the structural sheaf $\cO_{\widetilde{X}}$. This implies that the full transform $\sigma^*(\cI)$ is principal and generated by the sheaf of ideal of a divisor whose components are the exceptional divisors.
The actual process of desingularization is controlled by some
invariant and is often achieved before $(X,{\cI},E,1)$ has been resolved (Proposition \ref{pr: resol})
.

(2)$\Rightarrow$(3)
Let $Y\subset X$ be an
irreducible  subvariety. Assume there is a principalization of sheaves of ideals ${\cI}_Y$ subject to conditions (a) and (b) from Theorem \ref{th: 1}.
 Then in the course of the
principalization of ${\cI}_Y$ the strict transform $Y_i$
of $Y$ on some $X_i$ is the center of a blow-up. At this stage $Y_i$ is nonsingular and has simple normal crossing with the exceptional divisors. In the algorithm this moment is detected by some invariant.

(3)$\Rightarrow$(4) Every algebraic variety admits locally an embedding into an affine space.  Thus we can show that the existence of canonical embedded desingularization independent of the embedding defines a canonical desingularization. The patching of local desingularizations is controlled by an invariant independent of embeddings into smooth ambient varieties, provided the dimensions of the ambient varieties are the same.

\bigskip

\subsection{Functorial properties of multiple test  blow-ups} 
\begin{lemma}
Let $\phi: X' \to X$ be a smooth morphism. Then 
\begin{enumerate}
\item For any sheaf of ideals $\cI$ on $X$  and any $x'\in X'$, $x=\phi(x')\in X$ we have   $\ord_{x'}(\phi^*({\cI}))=\ord_{x}({\cI})$.
\item Let $E$ be a  set of divisors with smooth disjoint components such that all components of all divisors have simultaneously simple normal crossings. Then the inverse images $\phi^{-1}(D)$ of divisors  $D\in E$ have disjoint components and all the components of the divisors $\phi^{-1}(D)$ have simultaneously simple normal crossings.
\end{enumerate}\end{lemma}
\noindent{\bf Proof.} The assertions can be verified locally. 
Assume $\phi$ is of relative dimension $r$. Then it  factors (locally) as $\phi=\pi{\psi}$ where $\psi:U_{x'}\to X\times{\bf A}^r$  is \'etale  and $\pi:X\times{\bf A}^r  \rightarrow X$ is the natural projection. Let $x'':=\psi(x')$.
Since $\psi$ defines a formal analytic isomorphism $\psi^*:\widehat{\cO}_{x'',X\times{\bf A}^r}\simeq \widehat{\cO}_{x',X}$ the assertion of the lemma are satisfied  for $\psi$.  Moreover they are satisfied for the natural projection $\pi$ and for the composition $\phi=\pi\psi$.\qed

\begin{proposition} \label{pr: funct} Let $X_i$ be a multiple test  blow-up of a  marked ideal  $(X,{\cI},E,\mu)$ defining a sequence of marked ideals $(X_i,{\cI}_i,E_{i},\mu)$.
 Given a smooth morphism $\phi: X'\to X$, the induced sequence 
 $X'_i:=X'\times_{X}X_i$ is a multiple test  blow-up of $(X',{\cI}',E',\mu)$ such that  
\begin{enumerate}
 \item $\phi$ lifts to smooth morphisms $\phi_i: X'_i\to X_i$.
 
 \item $(X'_i)$ defines a sequence of marked ideals $(X'_i,{\cI}'_i,E'_{i},\mu)$ where
 ${\cI}'_i=\phi_i^*({\cI}_i)$,  the divisors in  $E'_i$ are the inverse images of the divisors in $E_i$ and the order on $E'_i$ is defined by the order on $E_i$.
\item If $(X_i)$ is a resolution of  $(X,{\cI},E,\mu)$ then  $(X'_i)$ is an extension of a resolution of $(X',{\cI}',E',\mu)$.
\end{enumerate}
\end{proposition}
\noindent{\bf Proof.} Induction on $i$. The pullback $\sigma'_{i+1}:X'_{i} \leftarrow X'_{i+1}$ of the blow-up $\sigma_{i+1}:X_{i} \leftarrow X_{i+1}$ by the smooth morphism $\phi_i: X'_i\to X_i$ is either the blow-up of the smooth center $C'_{i}=\sigma_i^{-1}(C_{i})$ or an isomorphism if $\sigma_i^{-1}(C)=\emptyset$. Since $\phi_i: X'_i\to X_i$ is smooth for any $x'\in X'_i$ and $x=\phi_i(x)$, $\ord_{x'}({\cI}'_i)=\ord_{x'}(\phi_i^*({\cI}_i))=\ord_{x}({\cI}_i)$.
Thus $C'_i\subset\supp(\cI'_i,\mu)$. Moreover $C'_i$ has simple normal crossings with $E'_i$. It is left to show that the transformation rules for the sheaves of ideals $\cI'_i$ and sets of divisors $E'_i$ are carried over by the induced smooth morphisms $\phi_{i}$. 
Note that the inverse image  of the exceptional divisor $D_{i+1}$ of $\sigma_{i+1}$ is the exceptional divisor $D'_{i+1}=\phi^{-1}_i(D_i)$ of $\sigma'_{i+1}$ . Thus we have
 \begin{eqnarray} &E'_{i+1}=(\sigma'_i)^c(E'_i)\cup \{D'_{i+1}\}=(\sigma'_i)^c(\phi_i^{-1}(E_i))\cup \phi_{i+1}^{-1}\{D_{i+1}\}=\phi_{i+1}^{-1}((\sigma_i)^c(E_i)\cup\{D_{i+1}\})=\phi_{i+1}^{-1}(E_{i+1})\nonumber \\
&{\cI}'_{i+1}=(\sigma'_{i+1})^c({\cI}'_i)=(\sigma'_{i+1})^*({\cI}'_i)\cI(D'_{i+1})^{-\mu}=(\sigma'_{i+1})^*(\phi_i^*({\cI}_i))\cdot \phi_{i+1}^*((\cI(D_{i+1}))^{-\mu})=\nonumber \\ &\phi_{i+1}^*(\sigma_{i+1}^*({\cI}_i)\cdot (\cI(D_{i+1}))^{-\mu})=
\phi_{i+1}^*(\sigma_{i+1}^c({\cI}_i)=\phi_{i+1}^*({\cI}_{i+1})\nonumber
\end{eqnarray}  \qed

\begin{definition}
We say that the above multiple test  blow-up  $(X'_i)$ is {\it induced} by $\phi_i$ and $X$.
We shall denote $(X'_i)$ and the corresponding marked ideals $(X',{\cI}',E',\mu)$ by
$$\phi^*(X_i):=X'_i,\quad \phi^*(X_i,{\cI}_i,E_i,\mu):=(X'_i,{\cI'_i},{E'_i},\mu).$$ 
\end{definition}

\bigskip
\subsection{Equivalence relation for marked ideals}

Let us introduce the following equivalence relation for marked ideals:
\begin{definition} Let  $(X,{\cI},E_{\cI},\mu_{\cI})$ and $(X,{\cJ},E_{\cJ},\mu_{\cJ})$ be two marked ideals 
on  the smooth variety $X$. Then

\centerline{$(X,{\cI},E_{\cI},\mu_{\cI})\simeq (X,{\cJ},E_{\cJ},\mu_{\cJ})$} 
\noindent if  
\begin{enumerate}
\item $E_{\cI}=E_{\cJ}$ and the orders on $E_{\cI}$ and on $E_{\cJ}$ coincide.
\item $\supp(X,{\cI},E_{\cI},\mu_{\cI})=\supp(X,{\cJ},E_{\cJ},\mu_{\cJ}).$
\item All the multiple test  blow-ups $ 
X_0=X\buildrel \sigma_1 \over\longleftarrow X_1 
\buildrel \sigma_2 \over\longleftarrow \ldots  \longleftarrow X_i 
\buildrel \sigma_{r} \over\longleftarrow\ldots \longleftarrow X_r$
of $(X,{\cI},E_{\cI},\mu_{\cI})$  are exactly the multiple test  blow-ups of $(X,{\cJ},E_{\cJ},\mu_{\cJ})$ and moreover
 we have
$$\supp(X_i,{\cI}_i,E_i,\mu_{\cI})=\supp(X_i,{\cJ}_i,E_i,\mu_{\cJ}).$$
\end{enumerate}
\end{definition}
\begin{example}
For any $k\in {\bf N}$, $({\cI},\mu)\simeq ({\cI}^k,k\mu)$.
\end{example}
\begin{remark} The marked ideals considered in this paper satisfy a stronger equivalence condition:
For any smooth morphisms $\phi: X' \to X$, $\phi^*(\cI,\mu)\simeq\phi^*(\cJ,\mu)$.
This condition will follow and is not added in the definition. 
 \end{remark}


\bigskip
\subsection{Ideals of derivatives}
Ideals of derivatives were first introduced and studied in the resolution context by Giraud.
Villamayor developed and applied this language to his {\it basic objects}.

\begin{definition}(Giraud, Villamayor) Let ${\cI}$ be a coherent sheaf of ideals on a smooth variety $X$. By the {\it  first derivative} (originally {\it extension}) ${\cD}({\cI})$ of $\cI$ we mean the coherent sheaf of ideals generated by all functions
$f\in {\cI}$ with their first derivatives. Then the {\it i-th derivative} ${\cD}^i({\cI})$ is defined to be ${\cD}({\cD}^{i-1}({\cI}))$. If $({\cI},\mu)$ is a marked ideal and $i\leq \mu$ then we define
$${\cD}^i({\cI},\mu):=({\cD}^i({\cI}),\mu-i).$$
\end{definition}
Recall that on a smooth variety $X$ there is a locally free sheaf of differentials $\Omega_{X/K}$ over $K$ generated locally by $du_1,\ldots, du_n$ for a set of local parameters $u_1,\ldots, u_n$. The dual sheaf of derivations $\Der_K(\cO_X)$ is locally generated by the derivations $\frac{\partial }{\partial u_i}$.
Immediately from the definition we observe that ${\cD}({\cI})$ is a coherent sheaf
defined locally by  generators $f_j$ of ${\cI}$ and all their partial derivatives $\frac{\partial f_j}{\partial u_i}$. We see by induction that ${\cD}^i({\cI})$ is a coherent sheaf
defined locally by the generators $f_j$ of ${\cI}$ and their derivatives $\frac{\partial^{|\alpha|} f_j}{\partial u^\alpha}$ for all multiindices $\alpha=(\alpha_1,\ldots,\alpha_n)$, where $|\alpha|:=\alpha_1+\ldots+\alpha_n\leq i$.

\begin{remark} In characteristic $p$ the partial derivatives $\frac{1}{\alpha!}\frac{\partial^{|\alpha|} }{\partial u^\alpha}$ (where $\alpha!=\alpha_1!\cdot\ldots\cdot\alpha_n!$) are well defined and are called the 
Hasse-Dieudonn\'e  derivatives. They should be used in the definition of the derivatives of marked ideals. One of the major sources of problems is that unlike in characteristic zero 
$${\cD}^i({\cD}^j({\cI}))\subsetneq {\cD}^{i+j}({\cI}).$$
\end{remark}

\begin{lemma}(Giraud, Villamayor) \label{le: Vi1}
For any  $i\leq\mu-1$,  $$\supp({\cI},\mu)=\supp({\cD}^i({\cI}),\mu-i)).$$
In particular \quad $\supp({\cI},\mu)=\supp({\cD}^{\mu-1}({\cI}),1)=V({\cD}^{\mu-1}({\cI})$\quad is a closed set.
 \qed
\end{lemma}

\bigskip
We write $({\cI},\mu)\subset ({\cJ},\mu)$ if $\cI\subset {\cJ}$.
\begin{lemma}(Giraud,Villamayor) \label{le: inclusions} Let $({\cI},\mu)$ be a marked ideal 
and $C\subset\supp ({\cI},\mu)$ be a smooth center and $r\leq \mu$. 
Let $\sigma: X\leftarrow
X'$ be a blow-up at $C$. Then $$\sigma^{\rm c}({\cD}^r({\cI},\mu)) \subseteq {\cD}^r(\sigma^{\rm c}({\cI},\mu)).$$
\end{lemma}
\noindent{\bf Proof.} 
First assume that $r=1$.
Let $u_1,\ldots,u_n$ denote the local parameters at $x$. Then the local
parameters at $x'\in\sigma^{-1}(x)$ are of the form 
$u'_i=\frac{u_i}{u_m}$ for $i<m$ and $u'_i=u_i$ for  $i\geq m$, where $u_m=u'_m=y$ denotes the local equation of the 
exceptional divisor.

The 
derivations $\frac{\partial{}}{\partial{u_i}}$ of $\cO_{x,X}$ extend to the derivations of the rational field  
$K(X)$. Note also that 
\begin{eqnarray}
&\frac{\partial{u'_j}}{\partial{u_i}}=\frac{\delta_{ij}}{u_m} , \quad i<m,  1\leq j \leq n; \quad \quad \frac{\partial{u'_j}}{\partial{u_m}}=-\frac{1}{u^2_m}u_j,\quad j< m; \quad\quad 
 \frac{\partial{u'_m}}{\partial{u_m}}=1,\quad\quad \frac{\partial{u'_j}}{\partial{u_m}}=0,  j> m; \nonumber\\
 &\frac{\partial{u'_i}}{\partial{u_j}}=\delta_{ij}, \quad i \geq m. \nonumber
\end{eqnarray} 

This gives 
\begin{eqnarray}
&\frac{\partial{}}{\partial{u_i}}=\frac{{1}}{u_m}\frac{\partial{}}{\partial{u'_i}}= \frac{{1}}{y}\frac{\partial{}}{\partial{u'_i}}, \quad 1\leq i < m; \quad\quad 
\frac{\partial{}}{\partial{u_m}}=-\frac{{1}}{y}(u'_1\frac{\partial{}}{\partial{u'_1}}+\ldots
+u'_{m-1} \frac{\partial{}}{\partial{u'_{m-1}}}-u'_m\frac{\partial{}}{\partial{u'_m}}),\nonumber\\
&\frac{\partial{}}{\partial{u'_i}}=\frac{\partial{}}{\partial{u_i}},   m<i \leq n. \nonumber
\end{eqnarray}
We see that  any derivation $D$ of $\cO_X$  induces a derivation $y\sigma^*(D)$ of $\cO_{X'}$.  Thus the sheaf $y\sigma^*(\\Der_K(\cO_X))$ of such derivations is a subsheaf of $\\Der_K(\cO_{X'})$ locally generated by \begin{equation}\frac{\partial{}}{\partial{u'_i}} , i<m ;  \quad y\frac{\partial{}}{\partial{y}}, \quad \mbox{and} \quad 
{y}\frac{\partial{}}{\partial{u'_i}}, i>m.\nonumber\end{equation}  
In particular
$y\sigma^*({\cD}({\cI}))\subset {\cD}(\sigma^*({\cI}))$. 
For any sheaf  of ideals ${\cJ}$ on $X'$ denote by $y\sigma^*({\cD})({\cJ})\subset {\cD}({\cJ})$
the ideal generated by ${\cJ}$ and  the derivatives  $D'(f)$, where $f\in \cJ$ and $D' \in y\sigma^*(\\Der_K(\cO_X))$. 
Note that for any $f\in \cJ$ and $D' \in y\sigma^*(\\Der_K(\cO_X))$,  $y$ divides $D'(y)$ and \begin{equation}D'(yf)=yD'(f)+D'(y)f\quad \in \quad y\sigma^*({\cD})(\cJ)+y\cJ=y\sigma^*({\cD})(\cJ).\nonumber \end{equation} Consequently, $y\sigma^*({\cD})(y{\cJ})\subseteq yy\sigma^*({\cD})({\cJ})$  and more generally $y\sigma^*({\cD})(y^\mu{\cJ})\subseteq y^{\mu}y\sigma^*({\cD})({\cJ})$. 
Then 
$$y\sigma^*({\cD}({\cI}))\subseteq y\sigma^*({\cD})(\sigma^*({\cI}))=y\sigma^*({\cD})(y^{\mu}\sigma^{\rm c}({\cI}))\subseteq y^{\mu}y\sigma^*({\cD})(\sigma^{\rm c}({\cI}))\subseteq y^{\mu}{\cD}(\sigma^{\rm c}({\cI})).$$
Then $$\sigma^{\rm c}({\cD}({\cI}))=y^{-\mu+1}\sigma^*({\cD}({\cI}))\subseteq {\cD}(\sigma^{\rm c}({\cI})).$$
Assume now that  $r$ is arbitrary. Then $C\subset\supp ({\cI},\mu)=\supp({\cD}^i({\cI},\mu))$ for $i\leq r$ and by induction on $r$, $$\sigma^{\rm c}({\cD}^{r}{\cI})=\sigma^{\rm c}({\cD}({\cD}^{r-1}({\cI})))\subseteq {\cD}(\sigma^{\rm c}{\cD}^{r-1}({\cI}))\subseteq {\cD}^{r}(\sigma^{\rm c}({\cI})).\qed$$

As a corollary from Lemma  \ref{le: inclusions} we prove the following

\begin{lemma} \label{le: inclusions2} A  multiple test  blow-up $(X_i)_{0\leq i\leq k}$ of $({\cI},\mu)$ is a multiple test  blow-up of ${\cD}^j({\cI},\mu)$ for $0\leq j\leq \mu$ and $$[{\cD}^j({\cI},\mu)]_k\subset {\cD}^j({\cI}_k,\mu).$$
\end{lemma}
\noindent {\bf Proof.} Induction on $k$. For $k=0$ evident. Let $\sigma_{k+1}:X_{k}\leftarrow X_{k+1}$ denote the blow-up with  a center 
$C_{k}\subseteq \supp({\cI}_k,\mu)=\supp({\cD}^j({\cI}_k,\mu))\subseteq \supp([{\cD}^j({\cI},\mu)]_k)$.
 Then by induction $[{\cD}^j({\cI},\mu)]_{k+1}=\sigma^{\rm c}_{k+1}
[{\cD}^j({\cI},\mu)]_k\subseteq \sigma^{\rm c}_{k+1}({\cD}^j({\cI}_k,\mu))$. Lemma \ref{le: inclusions} gives $\sigma^{\rm c}_{k+1}({\cD}^j({\cI}_k,\mu))\subseteq {\cD}^j\sigma_{k+1}^{\rm c}({\cI}_k,\mu)=
{\cD}^j({\cI}_{k+1},\mu)
$. \qed

\begin{lemma} \label{le: etale}  Let $\phi: X'\to X$ be an \'etale  morphism of smooth varieties   and ${\cI}$ be a sheaf of ideals on $X$. 

Then $$\phi^*({\cD}({\cI}))={\cD}(\phi^*({\cI})).$$
\end{lemma}
\noindent{\bf Proof.} 
Since $\phi$ is \'etale, for any points $x'\in X'$ and $\phi(x')=x$ we have  $\widehat{\cO}_{x',X'}\simeq\widehat{\cO}_{x,X}$, and 
$\widehat{\phi}^*(\widehat{{\cD}({\cI})})_{x}=({\cD}\widehat{(\phi^*({\cI})}))_{x'}$. Since $\widehat{\cO}_{x,X}$ is faithfully flat over ${\cO}_{x,X}$ and $\widehat{({\cD}({\cI}))}={{\cD}({\cI})}\cdot\widehat{\cO}_{x,X}$ we get  the equality of stalks $\phi^*({\cD}({\cI}))_x={\cD}(\phi^*({\cI}))_{x'}$ which determines the equality of sheaves $\phi^*({\cD}({\cI}))={\cD}(\phi^*({\cI}))$.
 \qed

\subsection{Hypersurfaces of maximal contact} 

The concept of the {\it hypersurfaces of maximal contact} is one of the key points of this proof. It was
originated by Hironaka, Abhyankhar and Giraud and developed in the papers of Bierstone-Milman and Villamayor.

In our terminology we are looking for a smooth hypersurface containing the supports of marked ideals and whose strict transforms under multiple test  blow-ups contain the supports of the induced marked ideals. Existence of such hypersurfaces allows a reduction of the resolution problem to  codimension 1.

First we introduce marked ideals which locally admit  hypersurfaces of maximal contact.

\begin{definition}(Villamayor (see \Vi))
We say that a marked ideal $({\cI},\mu)$  is of {\it maximal order} (originally {\it simple basic object}) if $\max\{\ord_x({\cI})\mid x\in X\}\leq \mu$ or equivalently ${\cD}^\mu({\cI})=\cO_X$.
\end{definition}

\begin{lemma}(Villamayor (see \Vi)) Let $({\cI},\mu)$ be a marked ideal of maximal order
and $C\subset\supp ({\cI},\mu)$ be a smooth center. Let $\sigma: X\leftarrow
X'$ be a blow-up at $C\subset\supp({\cI},\mu)$. Then $\sigma^{\rm c}({\cI},\mu)$ is of maximal order. 
\end{lemma}
\noindent {\bf Proof.} If $({\cI},\mu)$ is a marked ideal of maximal order then ${\cD}^\mu({\cI})=\cO_X$.
Then by  Lemma \ref{le: inclusions},  ${\cD}^\mu(\sigma^{\rm c}({\cI},\mu))\supset \sigma^{\rm c}({\cD}^\mu({\cI}),0)=\cO_X$. \qed

\begin{lemma}(Villamayor (see \Vi)) If $({\cI},\mu)$ is a marked ideal of maximal order and $0\leq i \leq \mu$ then ${\cD}^{i}({\cI},\mu)$
is of maximal order. 
\end{lemma}
\noindent {\bf Proof.} ${\cD}^{\mu-i}({\cD}^{i}({\cI},\mu))={\cD}^{\mu}({\cI},\mu)=\cO_X$.\qed

\begin{lemma}(Giraud)\label{le: Gi} Let $({\cI},\mu)$ be the marked ideal of maximal order.  Let $\sigma: X\leftarrow
X'$ be a blow-up at a smooth center $C\subsetneq\supp({\cI},\mu)$. Let $u\in {\cD}^{\mu-1}({\cI},\mu)(U)$ be a function of multiplicity one on $U$, that is, for any $x\in V(u)$, $\ord_x(u)=1$. In particular $\supp({\cI},\mu)\cap U \subset V( u)$. Let $U'\subset \sigma^{-1}(U)\subset X'$ be an open
set where the exceptional divisor is described by $y$.  Let $u':=\sigma^{\rm c}(u)=y^{-1}\sigma^*(u)$ be the controlled transform of $u$.
Then 
\begin{enumerate}
\item $u'\in {\cD}^{\mu-1}(\sigma^{\rm c}({\cI}_{|U'},\mu)).$ 
\item $u'$ is a function of multiplicity one on $U'$. 
\item $V(u')$ is the restriction of the strict transform of $V(u)$ to $U'$.

\end{enumerate}
\end{lemma}
\noindent {\bf Proof.} (1) $u'=\sigma^{\rm c}(u)=u/y\in\sigma^{\rm c}({\cD}^{\mu-1}({\cI}))\subset {\cD}^{\mu-1}(\sigma^{\rm c}({\cI}))$. 

(2) Since $u$ was one of the local parameters describing the center of blow-ups, $u'=u/y$ is a parameter, that
is, a function of order one.

 (3) follows from (2). \qed
\begin{definition} We shall call a function $$u\in T({\cI})(U):={\cD}^{\mu-1}({\cI}(U))$$ \noindent of multiplicity one a {\it tangent direction} of $({\cI},\mu)$ on $U$.
\end{definition}
As a corollary from the above we obtain the following lemma:
\begin{lemma}(Giraud)\label{le: Giraud} 
Let $u\in T({\cI})(U)$ be a tangent direction of $({\cI},\mu)$ on $U$. Then  for any multiple test  blow-up  $(U_i)$ of $({\cI}_{|U},\mu)$ all the supports  of the induced marked ideals  $\supp(\cI_i,\mu)$ are contained in the strict transforms $V(u)_i$ of $V(u)$.  \qed
\end{lemma} 
\begin{remarks}
\begin{enumerate}
\item Tangent directions are functions defining locally hypersurfaces of maximal contact.
\item  The main problem leading to complexity of the proofs is that of noncanonical choice of
the tangent directions. We overcome this difficulty by introducing {\it homogenized ideals}.
\end{enumerate}
\end{remarks}

\begin{lemma}(Villamayor) \label{le: codim1} Let $({\cI},\mu)$ be the marked ideal of maximal order whose support is of codimension $1$. Then all  codimension one components of $\supp({\cI},\mu)$ are smooth and isolated.  After the blow-up $\sigma: X\leftarrow
X'$ at such a component $C\subset\supp({\cI},\mu)$ the induced support $\supp({\cI}',\mu)$ does not intersect  the  exceptional divisor of $\sigma$.

\end{lemma}
\noindent{\bf Proof.} By the previous lemma there is a tangent direction $u\in {\cD}^{\mu-1}({\cI})$ whose zero set is smooth and contains $\supp({\cI},\mu)$. Then ${\cD}^{\mu-1}({\cI})=(u)$ and ${\cI}$ is locally described as
${\cI}=(u^\mu)$. The blow-up at the component $C$ locally defined by $u$ transforms $({\cI},\mu)$ to $({\cI}',\mu)$, where 
$\sigma^*({\cI})=y^\mu\cO_{X}$, and
${\cI'}=\sigma^{\rm c}(\cI)=y^{-\mu}\sigma^*({\cI})=\cO_{X}$, where $y=u$ describes the exceptional divisor. \qed 

\begin{remark} Note that the blow-up of codimension one components  is an isomorphism. However it defines a nontrivial transformation of marked ideals. In the actual desingularization process this kind of blow-up may occur for some marked ideals induced on subvarieties of ambient varieties. Though they define isomorphisms of those subvarieties they determine blow-ups of ambient
varieties which are not isomorphisms.
\end{remark}

\subsection{Arithmetical operations on marked ideals}
In this sections all marked ideals are defined for the smooth variety $X$ and the same
set of exceptional divisors $E$.
Define the following operations of addition and multiplication of marked ideals:
\begin{enumerate}
\item $({\cI},\mu_{\cI})+({\cJ},\mu_{\cJ}):=({\cI}^{\mu_{\cI}}+{\cJ}^{\mu_{\cI}},\mu_{\cI}\mu_{\cJ}),$
\noindent or more generally (the operation of addition is not associative)
$$({\cI}_1,\mu_1)+\ldots+({\cI}_m,\mu_m):=({\cI}_1^{\mu_2\cdot\ldots\cdot\mu_m}+
{\cI}_2^{\mu_1\mu_3\cdot\ldots\cdot\mu_m}+\ldots+{\cI}_m^{\mu_1\ldots\mu_{k-1}} ,\mu_1\mu_2\ldots\mu_m).$$
\item $(\cI,\mu_\cI) \cdot (\cJ,\mu_\cJ):=(\cI\cdot J,\mu_\cI+\mu_\cJ)$.
\end{enumerate}
\begin{lemma} \label{le: operations}
\begin{enumerate}
\item 
$\supp(({\cI}_1,\mu_1)+\ldots+({\cI}_m,\mu_m))= \supp({\cI}_1,\mu_1)\cap\ldots\cap\supp({\cI}_m,\mu_m)$. Moreover  multiple test  blow-ups $(X_k)$ of  $({\cI}_1,\mu_1)+\ldots+({\cI}_m,\mu_m)$ are exactly those which are simultaneous multiple test  blow-ups for all $({\cI}_j,\mu_j)$  and for any $k$ we have the equality for the controlled transforms $({\cI}_j,\mu_\cI)_k$ $$({\cI}_1,\mu_1)_k+\ldots+({\cI}_m,\mu_m)_k=[({\cI}_1,\mu_1)+\ldots+({\cI}_m,\mu_m)]_k$$
\item
$$\supp(\cI,\mu_\cI)\cap\supp(\cJ,\mu_\cJ)\supseteq\supp((\cI,\mu_\cI)\cdot (\cJ,\mu_\cJ)).$$
Moreover any simultaneous multiple test  blow-up $X_i$ of  both ideals $(\cI,\mu_\cI)$ and $(\cJ,\mu_\cJ)$ is a multiple test  blow-up for 
$(\cI,\mu_\cI)\cdot (\cJ,\mu_\cJ)$, and for the controlled transforms $(\cI_k,\mu_\cI)$ and $(\cJ_k,\mu_\cJ)$ we have the equality 
$$(\cI_k,\mu_\cI)\cdot(\cJ_k,\mu_\cJ)=[(\cI,\mu_\cI)\cdot(\cJ,\mu_\cJ)]_k.$$
\end{enumerate}
\end{lemma}
\noindent{\bf Proof.}  

(1) To simplify notation we restrict ourselves to the case of two marked ideals. The proof for $n>2$ marked ideals is exactly the same. We have

$\supp(({\cI},\mu_{\cI})+({\cJ},\mu_{\cJ}))=\supp({\cI}^{\mu_{\cJ}}+{\cJ}^{\mu_{\cI}},\mu_{\cI}\mu_{\cJ})=\supp({\cI}^{\mu_{\cJ}},\mu_{\cI}\mu_{\cJ})
\cap\supp({\cJ}^{\mu_{\cI}},\mu_{\cI}\mu_{\cJ})$=\\
$\supp({\cI},\mu_{\cI})\cap\supp({\cJ},\mu_{\cJ})$.

Suppose now all multiple test  blow-ups of $({\cI},\mu_{\cI})+({\cJ},\mu_{\cJ})$ of length $k\geq 0$  are  exactly simultaneous multiple test  blow-ups for $({\cI},\mu_{\cI})$  and $({\cJ},\mu_{\cJ})$ and  $[({\cI},\mu_{\cI})+({\cJ},\mu_{\cJ})]_{k}=({\cI}_k,\mu_{\cI})+({\cJ}_k,\mu_{\cJ})$.
 Let $\sigma_{k+1}$ denote a blow-up of smooth center $C_k$ contained in $\supp[({\cI},\mu_{\cI})+({\cJ},\mu_{\cJ})]_k$. Then
\begin{eqnarray}
\noindent &[({\cI},\mu_{\cI})+({\cJ},\mu_{\cJ})]_{k+1}=\sigma_{k+1}^{\rm c}(({\cI}_k,\mu_{\cI})+({\cJ}_k,\mu_{\cJ}))=(y^{-\mu_{\cI}\mu_{\cJ}}\sigma_{k+1}^*({\cI}_k^{\mu_{\cJ}}+{\cJ}_k^{\mu_{\cI}},\mu_{\cI}\mu_{\cJ}))= \nonumber\\
\noindent & (y^{-\mu_{\cI}\mu_{\cJ}}\sigma_k^*({\cI}_k)^{\mu_{\cJ}},\mu_{\cI}\mu_{\cJ})+
(y^{-\mu_{\cI}\mu_{\cJ}}\sigma_k^*({\cJ}_k)^{\mu_{\cI}},\mu_{\cI}\mu_{\cJ})=
(y^{-\mu_{\cI}}\sigma_k^{*}({\cI}_k),{\mu_{\cI}})+(y^{-\mu_{\cJ}}\sigma_k^{*}({\cJ}_k),\mu_{\cJ})=\nonumber \\
\noindent &\sigma_k^{\rm c}({\cI}_k,{\mu_{\cI}})+\sigma_k^{\rm c}({\cJ}_k,\mu_{\cJ})=({\cI}_{k+1},{\mu_{\cI}})+({\cJ}_{k+1},\mu_{\cJ}).\nonumber
\end{eqnarray}

 (2)  If $\ord_x(\cI)\geq\mu_{\cI}$ and  $\ord_x(\cJ)\geq\mu_{\cJ}$ then $\ord_x(\cI\cdot\cJ)\geq\mu_{\cI}+\mu_\cJ$. This implies that \\ $\supp(\cI,\mu_\cI)\cap\supp(\cJ,\mu_\cJ)\supseteq\supp((\cI,\mu_\cI)\cdot (\cJ,\mu_\cJ)).$
Suppose now all simultaneous multiple test  blow-ups of $({\cI},\mu_{\cI})$  and $({\cJ},\mu_{\cJ})$ of length $k\geq 0$  are  multiple test  blow-ups for  $({\cI},\mu_{\cI})\cdot({\cJ},\mu_{\cJ})$ and 
 there is equality \begin{equation} \nonumber ({\cI}_k,\mu_{\cI})\cdot({\cJ}_k,\mu_{\cJ})=[({\cI},\mu_{\cI})\cdot({\cJ},\mu_{\cJ})]_{k}.
 \end{equation}
 Let $\sigma_{k+1}$ denote the blow-up of a smooth center $C_k$ contained in $\supp({\cI}_k,\mu_{\cI})\cdot({\cJ}_k,\mu_{\cJ})\subseteq \supp({\cI},_k\mu_{\cI})\cap\supp({\cJ}_k,\mu_{\cJ})$. Then

\begin{eqnarray}
\noindent &[({\cI},\mu_{\cI})\cdot({\cJ},\mu_{\cJ})]_{k+1}=\sigma_{k+1}^{\rm c}(({\cI}_k,\mu_{\cI})\cdot({\cJ}_k,\mu_{\cJ}))=(y^{-(\mu_{\cI}+\mu_{\cJ})}\sigma_{k+1}^*({\cI}_k^\cdot {\cJ}_k), \mu_{\cI}+\mu_{\cJ}))\supseteq  \nonumber\\
\noindent  & (y^{-\mu_{\cI}}\sigma_k^*({\cI}_k),\mu_{\cI})\cdot
(y^{-\mu_{\cJ}}\sigma_k^*({\cJ}_k), \mu_{\cJ})=
\sigma_k^{\rm c}({\cI}_k,{\mu_{\cI}})\cdot\sigma_k^{\rm c}({\cJ}_k,\mu_{\cJ})=\nonumber ({\cI}_{k+1},{\mu_{\cI}})\cdot({\cJ}_{k+1},\mu_{\cJ}).\nonumber
\end{eqnarray}

\qed 
\begin{remark} The operation of multiplication of marked ideals is associative while the operation of addition is not. However we have the following lemma.
\end{remark}

\begin{lemma}
$(({\cI}_1,\mu_1)+({\cI}_2,\mu_2))+({\cI}_3,\mu_3)\simeq ({\cI}_1,\mu_1)+(({\cI}_2,\mu_2)+({\cI}_3,\mu_3))$

\end{lemma}
\noindent{\bf Proof.} It follows from Lemma \ref{le: operations} that the supports of the two marked ideals are the same. Moreover by the same lemma the supports remain the same after 
consecutive blow-ups of multiple test  blow-ups.
\qed

\bigskip
\subsection{Homogenized ideals and tangent directions}
Let $({\cI},\mu)$ be a marked ideal of maximal order. Set $T({\cI}):={\cD}^{\mu-1}{\cI}$.
By the {\it homogenized ideal} we mean $${\cH}({\cI},\mu):=({\cH}({\cI}),\mu)=({\cI}+{\cD}\cI\cdot T({\cI})+\ldots+{\cD}^i\cI\cdot T({\cI})^i+ \ldots+{\cD}^{\mu-1}\cI\cdot T({\cI})^{\mu-1},\mu)$$

\begin{lemma} \label{le: hom} Let $({\cI},\mu)$ be a marked ideal of maximal order.
\begin{enumerate}
 
 \item If $\mu=1$, then $({\cH}({\cI}),1)=(\cI,1)$.
\item ${\cH}({\cI})=({\cI}+{\cD}\cI\cdot T({\cI})+\ldots +{\cD}^i\cI\cdot T({\cI})^i+ \ldots +{\cD}^{\mu-1}\cI\cdot T({\cI})^{\mu-1}+{\cD}^{\mu}\cI\cdot T({\cI})^{\mu}+\ldots$.
\item $({\cH}({\cI}),\mu)=({\cI},\mu)+{\cD}(\cI,\mu)\cdot (T({\cI}),1)+\ldots +{\cD}^i(\cI,\mu)\cdot (T({\cI}),1)^i+ \ldots+ {\cD}^{\mu-1}(\cI,\mu)\cdot (T({\cI}),1)^{\mu-1}$.
\item If $\mu>1$ then ${\cH}(\cD({\cI},\mu))\supseteq \cD({\cH}({\cI},\mu))$.  \qed
\item ${T}(\cH({\cI},\mu))={T}({\cI},\mu)$.
\end{enumerate}
\end{lemma}

\begin{remark}
A homogenized ideal features two important properties:
\begin{enumerate}
\item It is equivalent to the given ideal. 
\item It ''looks the same'' from all possible tangent directions.
\end{enumerate}
By the first property we can   use the homogenized ideal to construct resolution via the Giraud Lemma \ref{le: Giraud}. By the second property such a construction does not depend on the choice of tangent directions.
\end{remark}

\begin{lemma} Let $({\cI},\mu)$ be a marked ideal of maximal order. Then
\begin{enumerate}
\item $({\cI},\mu)\simeq ({\cH}({\cI}),\mu)$.
\item For any multiple test  blow-up $(X_k)$ of $(\cI,\mu)$, \\$({\cH}({\cI}),\mu)_k=({\cI},\mu)_k+[{\cD}(\cI,\mu)]_k\cdot [(T({\cI}),1)]_k+\ldots [{\cD}^{\mu-1}(\cI,\mu)]_k\cdot +[(T({\cI}),1)]_k^{\mu-1}$.
\end{enumerate}
\end{lemma}

\noindent{\bf Proof.} Since $\cH(\cI)\supset \cI$, every multiple test  blow-up of $\cH(\cI,\mu)$ is a multiple test  blow-up of $(\cI,\mu)$. 
By Lemma \ref{le: inclusions2}, every multiple test  blow-up of  $(\cI,\mu)$ is a multiple test  blow-up for all ${\cD}^i(\cI,\mu)$ and consequently, by Lemma \ref{le: operations} 
 it is a simultaneous resolution of all 
$({\cD}^i(\cI)\cdot T({\cI})^i,\mu)=({\cD}^i(\cI),\mu-i)\cdot (T({\cI})^i,i)$ and \begin{eqnarray}
&\supp(\cH(\cI,\mu)_k)=\bigcap_{i=0}^{\mu-1} \supp({\cD}^i(\cI)\cdot T({\cI})^i,\mu)_k=\bigcap_{i=0}^{\mu-1} \supp({\cD}^i(\cI),\mu-i)_k\cdot (T({\cI})^i,i)_k\supseteq \nonumber \\
 &\bigcap_{i=0}^{\mu-1}\supp({\cD}^i(\cI,\mu))_k=\supp(\cI_k,\mu). \nonumber
 \end{eqnarray}
Therefore every multiple test  blow-up of  $(\cI,\mu)$ is a multiple test  blow-up of $\cH(\cI,\mu)$  and
by Lemmas \ref{le: hom}(3) and  \ref{le: operations}  we get (2).
\qed
\begin{lemma} \label{le: etale2} Let $\phi: X'\to X$ be a smooth morphism of smooth varieties 
and let $(X,{\cI},\emptyset,\mu)$ be a marked ideal. Then 
$$ \phi^*({\cH}({\cI}))={\cH}(\phi^*({\cI})).$$
\end{lemma}
\noindent{\bf Proof.} A direct consequence of Lemma \ref{le: etale} \qed
.

Although the following Lemmas  \ref{le: homo0} and \ref{le: homo} are used in this paper only in the case $E=\emptyset$ we formulate them in slightly more general versions.
\begin{lemma} \label{le: homo0} Let $(X,{\cI},E,\mu)$ be a marked ideal of maximal order. 
Assume there exist   tangent directions $u,v\in T(\cI,\mu)_x={\cD}^{\mu-1}({\cI},\mu)_x$ at $x\in\supp({\cI},\mu)$ which are transversal to $E$.
Then there exists an automorphism $\widehat{\phi}_{uv}$ of $\widehat{X}_x:=\Spec(\widehat{\cO}_{x,X})$ such that 
\begin{enumerate}
\item  $\widehat{\phi}_{uv}^*({\cH}\widehat{\cI})_x=({\cH}\widehat{\cI})_x.$ 
\item  $\widehat{\phi}_{uv}^*(E)=E.$
\item  $\widehat{\phi}_{uv}^*(u)=v.$
\item  $\supp(\widehat{\cI},\mu):=V(T(\widehat{\cI},\mu))$ is contained in the fixed point set of $\phi$. 
\end{enumerate}  
\end{lemma}
\noindent{\bf Proof.} (0) {\bf Construction of the automorphism $\widehat{\phi}_{uv}$.}

Find parameters $u_2,\ldots, u_n$ transversal to $u$ and $v$ such that
 $u=u_1,u_2,\ldots, u_n$  and $v,u_2,\ldots, u_n$ form  two sets of parameters at $x$ and divisors in $E$ are described by some parameters $u_i$ where $i\geq 2$.
Set \begin{equation}\widehat{\phi}_{uv}(u_1)=v, \quad \widehat{\phi}_{uv}(u_i)=u_i \quad\mbox{for} \quad i>1.\nonumber \end{equation}
(1)  Let $h:=v-u\in {T({\cI})}$. For any  $f\in \widehat{\cI}$, 
 \begin{equation} \widehat{\phi}_{uv}^*(f)=f(u_1+h,u_2,\ldots,u_n)= f(u_1,\ldots,u_n)+\frac{\partial{f}}{\partial{u_1}}\cdot h+ \frac{1}{2!}
 \frac{\partial^2{f}}{\partial{u_1^2}}\cdot h^2+\ldots +\frac{1}{i!}
 \frac{\partial^i{f}}{\partial{u_1^i}}\cdot h^i +\ldots\nonumber \end{equation}
 The latter element belongs to \begin{equation}\widehat{\cI}+{\cD}\widehat{\cI}\cdot \widehat{T({\cI})}+\ldots +{\cD}^i\widehat{\cI}\cdot \widehat{T({\cI})}^i+ \ldots +{\cD}^{\mu-1}\widehat{\cI}\cdot \widehat{T({\cI})}^{\mu-1}={\cH}\widehat{\cI}.\nonumber\end{equation}

Hence $\widehat{\phi}_{uv}^*(\widehat{\cI})\subset {\cH}\widehat{\cI}$. Analogously  $\widehat{\phi}_{uv}^*({\cD}^i\widehat{\cI})\subset
{\cD}^i\widehat{\cI}+{\cD}^{i+1}\widehat{\cI}\cdot T({\cI})+ \ldots +
{\cD}^{\mu-1}\widehat{\cI}\cdot \widehat{T({\cI})}^{\mu-i-1}={\cH}\widehat{{\cD}^iI}$.
In particular by Lemma \ref{le: hom}, $\widehat{\phi}_{uv}^*(\widehat{T({\cI})},1)\subset \cH(\widehat{T({\cI})},1)=(\widehat{T({\cI})},1) $. This gives

$$\widehat{\phi}_{uv}^*({\cD}^i\widehat{\cI}\cdot \widehat{T({\cI})}^i )\subset {\cD}^i\widehat{\cI}\cdot \widehat{T({\cI})}^i+ \ldots +{\cD}^{\mu-1}\widehat{\cI}\cdot \widehat{T({\cI})}^{\mu-1}\subset {\cH}\widehat{\cI}.$$
By the above  $\widehat{\phi}_{uv}^*({\cH}\widehat{\cI})_x\subset ({\cH}\widehat{\cI})_x$ and since the scheme is noetherian, $\widehat{\phi}_{uv}^*({\cH}\widehat{\cI})_x=({\cH}\widehat{\cI})_x.$

\bigskip
(2)(3) Follow from the construction.

(4) The fixed point set of $\widehat{\phi}^*_{uv}$ is defined by $u_i=\widehat{\phi}^*_{uv}(u_i)$, $i=1,\ldots,n$,  that is, $h=0$. But $h\in {\cD}^{\mu-1}({\cI})$ is $0$ on $\supp({\cI},\mu)$. \qed

\begin{lemma} \label{le: homo} ({\bf Glueing Lemma}) Let $(X,{\cI},E,\mu)$ be a marked ideal of maximal order for which there exist
tangent directions $u,v\in T({\cI},\mu)$ at  $x\in\supp({\cI},\mu)$ which are transversal to $E$.
Then there exist \'etale neighborhoods $\phi_{u},\phi_v: \overline{X}\to X$ of $x=\phi_u(\overline{x})=\phi_v(\overline{x}) \in X$,  where $\overline{x}\in \overline{X}$, such that 
\begin{enumerate}
\item  $\phi_{u}^*(X,{\cH}({\cI}),E,\mu)=\phi_{v}^*(X,{\cH}({\cI}),E,\mu)$. 
\item  $\phi_{u}^*(u)=\phi_{v}^*(v)$.

Set $(\overline{X},\overline{\cI},\overline{E},\mu):=\phi_{u}^*(X,{\cH}({\cI}),E,\mu)=\phi_{v}^*(X,{\cH}({\cI},E,\mu)).$
\item  For any $\overline{y}\in \supp(\overline{X},\overline{\cI},\overline{E},\mu)$, $\phi_u(\overline{y})=\phi_v(\overline{y})$.
\item  Let $(X_i)$ be a multiple test  blow-up of $({X},{\cI},\emptyset,\mu)$. Then
\begin{enumerate} 
\item The induced multiple test  blow-ups $\phi_u^*(X_i)$ and $\phi_v^*(X_i)$ of $(\overline{X},\overline{\cI},\overline{E},\mu)$ are the same (defined by the same centers).  Set $(\overline{X}_i):=\phi_u^*(X_i)=\phi_v^*(X_i)$ and let $\phi_{u i},\phi_{v i}:\overline{X}_i\to X_i$ be  the induced morphisms. Then
$\phi_{u i}^*(X_i,{\cH}({\cI})_i,E_i,\mu)=\phi_{v i}^*(X_i,{\cH}({\cI})_i,E_i,\mu)$.
\item Let $V(u)$ and $V(v)$ denote the hypersurfaces of maximal contact on $X$  and $V(u)_i$ and $V(v)_i$ be their strict transforms. Then
$\phi_{ui}^{-1}(V(u)_i)=\phi_{vi}^{-1}(V(v)_i)$.

\item For any  $\overline{y}_i\in \supp(\overline{X}_i,\overline{\cI}_i,\overline{E}_i,\mu)$  $\phi_{u i}(\overline{y}_i)=\phi_{v i}(\overline{y}_i).$
\end{enumerate}
\end{enumerate}
\end{lemma}

\noindent {\bf Proof.} 
(0) {\bf Construction of  \'etale neighborhoods  ${\phi}_{u}, {\phi}_{v}: U\to X$.} 

Let $U\subset X$ be an open subset for which there exist $u_2,\ldots, u_n$ which are transversal to $u$ and $v$ on $U$ such  that
 $u=u_1,u_2,\ldots, u_n$  and $v,u_2,\ldots, u_n$ form  two sets of parameters on $U$ and divisors in $E$ are described by some $u_i$, where $i\geq 2$.
 Let ${\bf A}^n$ be the affine space with coordinates $x_1,\ldots,x_n$.
 Construct first \'etale morphisms $\phi_1,\phi_2: U\to {\bf A}^n$ with  
 \begin{equation} {\phi}^*_{1}(x_i)=u_i \quad \textrm{for all } i \quad\quad \mbox{and} \quad
{\phi}^*_{2}(x_1)=v, \quad {\phi}^*_{2}(x_i)=u_i\quad\mbox{for} \quad i>1.\nonumber \end{equation}

\noindent  Consider the fiber product $U\times_{{\bf A}^n}U$ for  the morphisms $\phi_1$ and $\phi_2$. Let $\phi_u$, $\phi_v$ are the natural projections $\phi_u, \phi_v: U\times_{{\bf A}^n}U\to U$ such that $\phi_1\phi_u=\phi_2\phi_v$. Then define $\overline{X}$  to be an irreducible component of $U\times_{{\bf A}^n}U$ whose images $\phi_u(U)$ and $\phi_v(U)$  contain $x$ . Set
\begin{eqnarray} w_1:=\phi_u^*(u)=(\phi_1\phi_u)^*(x_1)=(\phi_2\phi_v)^*(x_1)=\phi_v^*(v), \nonumber\\
 w_i=\phi_u^*(u_i)=\phi_v^*(u_i)\quad  \textrm{for $i\geq 2$}. \nonumber \end{eqnarray}

(1) Let $h:=v-u$. By the above the morphisms $\phi_u$ and $\phi_v$ coincide on $\phi_u^{-1}(V(h))=\phi_v^{-1}(V(h))$.  If $\overline{y}\in \overline{X}$ be a point such that $\phi_u(\overline{y})\not\in\supp(X,{\cI},E,\mu)$ then $\overline{y}\not\in \supp(\phi_u^*(\cH(\cI))=\supp(\cH(\phi_u^*(\cI))$ and we have the equality of stalks $\phi_u^*(\cH(\cI)_{\overline{y}}=\cH(\phi_u^*(\cI)_{\overline{y}}=\cO_{\overline{X},\overline{y}}$. On the other hand $\phi_v(\overline{y})\not\in\supp(X,{\cI},E,\mu)$ and $\phi_v^*(\cH(\cI)_{\overline{y}}=\cH(\phi_v^*(\cI)_{\overline{y}}=\cO_{\overline{X},\overline{y}}$. 

Let $\overline{y}\in \overline{X}$ be a point such that $\phi_u(\overline{y})=y\in\supp(X,{\cI},E,\mu)$. Then $\phi_v(\overline{y})=y$. Denote by $(\widehat{\phi}_{v})_{\overline{y}}$ and $(\widehat{\phi}_{u})_{\overline{y}}$ the induced morphisms of the completions $\widehat{\overline{X}}_{\overline{y}}\to {\widehat{X}}_{{y}}$. We have the equality  $\widehat{\phi}_{uv}=(\widehat{\phi}_{v})_{\overline{y}}(\widehat{\phi}^{-1}_{u})_{\overline{y}}$ where
$\widehat{\phi}_{uv}$ is given as in the proof of Lemma \ref{le: homo0}.
 Consequently, for any such ${\overline{y}}$,
 \begin{equation}({\cH}\widehat{\cI})_{\overline{y}}=\widehat{\phi}^*_{uv}({\cH}\widehat{\cI})_{\overline{y}}=
 (\widehat{\phi}_{u})_{\overline{y}}(\widehat{\phi}^{-1}_{v})^*_{\overline{y}}({\cH}\widehat{\cI})_y.\nonumber\end{equation}
\begin{equation}{\phi}_{u}^*(\cH{\cI})_{\overline{y}}=(\widehat{\phi}_{u})^*_{\overline{y}}({\cH}\widehat{\cI})_{\overline{y}}=(\widehat{\phi}^{*}_{v})_{\overline{y}}({\cH}\widehat{\cI})_{\overline{y}}={\phi}_{v}^*(\cH{\cI})_{\overline{y}}.\nonumber \end{equation}
We get the equality of stalks  ${\phi}_{u}^*(\cH{\cI})_{\overline{y}}={\phi}_{v}^*(\cH{\cI})_{\overline{y}}$ for all points $y\in \overline{X}$ and for sheaves ${\phi}_{u}^*(\cH({\cI}))={\phi}_{v}^*(\cH({\cI}))$. Also ${\phi}_{u}^*(E)={\phi}_{v}^*(E)$ by construction and ${\phi}_{u}^*(T({\cI}))={\phi}_{v}^*(T({\cI}))$.

(2) Follows from the construction.

(3) Let $h:=v-u$.
The subset of $\overline{X}$ for which $\phi_1(x)=\phi_2(x)$ is described by $h=0$.  Consequently $\phi_u=\phi_v$ over $V(\phi_u^*(h))$. In particular these morphisms are equal over  $\supp({\cI},\mu)=\supp(\cH({\cI}),\mu)$.

(4)  Let $(X_i)$ be a multiple test blow-up of $(X,{\cI},E,\mu)$. Let $C_0\subset \supp(X,{\cI},E,\mu)$ be the center of $\sigma_1$. By (3), $\phi_u=\phi_v$ over $\supp(\cI,\mu)$. Fix a point $\overline{y}\in \supp(\phi_u^*(\cH(\cI)),\mu)$ and let
$y=\phi_u(\overline{y})=\phi_v(\overline{y})\in \supp(X,{\cI},E,\mu).$ 
Find parameters $u'_1=u_1,u'_2,\ldots,u'_n$ on  an affine neighborhood $U'$ of  $y$  such that 
divisors in $E$ are described by some $u_i$ for $i\geq 2$ and
$C_0$ is described  by $u'_1=u'_2=\ldots=u'_m=0$ for some $m\geq 0$.

Let $\overline{U}\subset \phi_u^{-1}(U)\cap\phi_v^{-1}(U)\subset \overline{X}$ be an affine neighbourhood of $\overline{y}$.

Let $\overline{J}$ be the ideal of $K[\overline{U}]$ generated by all functions $\phi^*_u(f)-\phi^*_v(f)$, where $f\in K[U']$. Then $(\phi_u^*(h))\subset \overline{J}$. On the other hand by definition for any point $\overline{z}\in V(\phi_u^*(h))\subset V(\overline{J})$ we have the equalities of the completions of stalks of the ideals
$$\widehat{\overline{J}}_{\overline{z},\overline{U}}=(\phi_u^*(u'_i)-\phi_v^*(u'_i))_{i=1,\ldots,n}=(\phi_u^*(h))\cdot\widehat{\cO}_{\overline{z},\overline{U}},$$ \noindent which implies the equalites of stalks of the ideals
$$(\phi_u^*(h))_{\overline{z},\overline{U}}= {\overline{J}}_{\overline{z},\overline{U}}.$$
 Finally $\overline{J}=(\phi_u^*(h))$ and consequently, for any $i=1,\ldots n$, we have 
$$\phi^*_u(u'_i)-\phi^*_v(u'_i)\in (\phi_u^*(h))\subset {\phi}_{u}^*(T({\cI}))(\overline{U})={\phi}_{v}^*(T({\cI}))(\overline{U}).$$

For simplicity denote $u'_i$ simply by $u_i$ and $U'$ by $U$.
Let $\overline{\sigma}_1:\overline{X}_1\to \overline{X}$ be the blow-up of $\overline{X}$ at $\overline{C}_0:=\phi_u^{-1}(C_0)=\phi_v^{-1}(C_0)\subset \overline{X}$. Then both morphisms $\phi_u$ and $\phi_v$ lift to the \'etale morphisms $\phi_{u 1},\phi_{v 1}: \overline{X}_1\to X_1$ by the universal property of a blow-up. 
Observe that $$\phi^*_{u 1}(T(\cI)_1,1)=\phi^*_{u 1}(\cI(D)^{-1}\cdot\sigma^*(T(\cI)))=\phi^*_{v 1}(\cI(D)^{-1}\cdot\sigma^*(T(\cI)))=\phi^*_{v 1}(T(\cI)_1).$$\noindent where $D$ denotes the exceptional divisor of $\sigma_1$.

Fix a point $\overline{y}_{1}\in \supp(\phi_{u 1}^*(\overline{\cI}_1),\mu)\subset \supp(\phi_{u 1}^*(T(\cI)_1),1)$ and $y_1=\phi_{u 1}(\overline{y}_{1})\in \supp(\cI_1,\mu)\subset \supp(T({\cI}_1),1)$ such that 
 $\overline{y}=\overline{\sigma}_1(\overline{y}_{1})\in \supp(\phi_u^*(\cI),\mu)$ and 
 $y=\sigma_1(y_1)\in \supp(\cI,\mu)$.
 Find parameters $u_1,u_2,\ldots,u_n$ at $y$, by replacing  $u_{2},\ldots,u_n$ if necessary by their linear combinations, 
such that
\begin{enumerate}
\item The  parameters at $y_1$ are given by $u_{i1}:=\frac{u_i}{u_m}$ for $1\leq i <m$ and $u_{i1}:=u_i$ for $i\geq m$ 
\item All divisors in $E_1$ through $y_1$ are defined by some $u_{i1}$. 
\end{enumerate}

Then $w_i:=\phi_u^*(u_i)$ for $i=1,\ldots,n$ define parameters at a point $\overline{y}$ such that the  parameters at 
$\overline{y}_{1}$ are given by $w_{i1}:=\frac{w_i}{w_m}$ for $i< m$ and $w'_{i1}:=w_i$ for $i\geq m$.

Let $U_{m1}\subset X_1$ be the neighbourhood defined by the parameter $u_m$. The subset
$U_{m1}\ni y_1$ is described by all points $z$ for which $(u_m/y_D)(z)\neq 0$, where $y_D$ is a local equation of the exceptional divisor of $\sigma_1$. Then a point $\overline{z}$ is in $\supp(\phi_{u 1}^*(T(\cI_1)),\mu)\cap \phi_u^{-1}(U_{m1})$ iff $(\phi_{u1}^*(u_m)/y_{\overline{D}})(\overline{z})\neq 0$,where $y_{\overline{D}}$ is a local equation of the exceptional divisor of $\overline{\sigma}_1$. But  for any $\overline{z}\in\supp(\phi_{u 1}^*(T(\cI_1)),1)$,  $$(\phi_{v}^*(u_m)/y_{\overline{D}})(x)=(\phi_{u}^*(u_m)/y_{\overline{D}})(x)+\phi_{u}^*(h)/y_{\overline{D}}(x)=(\phi_{u}^*(u_m)/y_{\overline{D}})(x),$$\noindent  where $h\in (T(\cI))(U)$ and $\phi_{u}^*(h)/y_{\overline{D}}\in \phi^*_{u 1}(T(\cI)_1)(\phi_{u1}^{-1}(U_m))$. This implies that
$$\supp(\phi_{u 1}^*(T(\cI_1)),\mu)\cap \phi_{u 1}^{-1}(U_{m1})=\supp(\phi_{v 1}^*(T(\cI_1)),\mu)\cap \phi_{v 1}^{-1}(U_{m1}).$$ Then the exceptional divisor of $\overline{\sigma}_1$ is described on  $\overline{U}_{m1}:=\phi_{u 1}^{-1}(U_{m1})\cap \phi_{v 1}^{-1}(U_{m1})$ and by ${\phi^*_{u 1}(u_{m})}$ and by ${\phi^*_{v 1}(u_{m})}$. 
\begin{eqnarray}\phi^*_{u 1}(u_{i 1})-\phi^*_{v 1}(u_{i 1})=
\phi^*_{u 1}(\frac{u_{i}}{u_m})-\phi^*_{u 1}(\frac{u_{i}}{u_m})=
\frac{\phi^*_{u 1}(u_{i})}{\phi^*_{u 1}(u_{m})}-\frac{\phi^*_{u 1}(u_{i})}{\phi^*_{v 1}(u_{m})}+
\frac{\phi^*_{u 1}(u_{i})}{\phi^*_{v 1}(u_{m})}-\frac{\phi^*_{u 1}(u_{i})}{\phi^*_{v 1}(u_{m})}=\nonumber\\
\frac{\phi^*_{u 1}(u_{i})}{\phi^*_{u 1}(u_{m})}\frac{(\phi^*_{v 1}(u_{m})-\phi^*_{u 1}(u_{m}))}{\phi^*_{v 1}(u_{m})}+
\frac{\phi^*_{u 1}(u_{i})-\phi^*_{v 1}(u_{i})}{\phi^*_{v 1}(u_{m})} \quad
\in \quad {\phi}_{u 1}^*(T({\cI})_1(\overline{U}_{m1})={\phi}_{v 1}^*(T({\cI}))_1(\overline{U}_{m1}).\quad \nonumber\end{eqnarray}
In particular  $u_{i 1}(\phi_{v 1}(\overline{y}_1))=\phi^*_{v 1}(u_{i 1})(\overline{y})=\phi^*_{u 1}(u_{i 1})(\overline{y})=0$. Thus $\phi_{v 1}(\overline{y}_1)=y_1=\phi_{u 1}(\overline{y}_1)$ is the only point in $\sigma_1^{-1}(y)$ for which all $u_{i1}$ are zero. This shows that  $\phi_{u 1}=\phi_{v 1}$ over $\supp(T({\cI})_1(U_{m1}),1)\supset \supp(\cI_1(U_{m1}),\mu)$ and thus over $\supp(T({\cI})_1\supset \supp(\cI_1,\mu)$. 

Set $v_{11}:=\frac{v}{u_m}$. Note that the functions $$\phi_{u1}^*(u_{11})=\phi_{u}^*(u)/\phi_{u}^*(u_m)=\phi_{v}^*(v)/\phi_{u}^*(u_m)\sim \phi_{v1}^*(v_{11})=\phi_{v}^*(v)/\phi_{v}^*(u_m)$$ \noindent are proportional (up  to an invertible function on $\overline{U}_{m1}$).
Also for any parameters $u_{i1}$ defining divisors in $E_1$ through $y_1$, we have
$$\phi_{u1}^*(u_{i1})=\phi_{u}^*(u_i)/\phi_{u}^*(u_m)\sim\phi_{v1}^*(u_{i1})/\phi_{u}^*(u_m).$$

Denote $U_{m1}$ and $\overline{U}_{m1}$ simply by $U_1$ and $\overline{U}_{1}$.  By induction on $k$ we show that $\phi_{u k}=\phi_{v k}$ over $\supp(T({\cI})_k,1)\supset \supp(\cI_k,\mu)$ and for any point $\overline{y}_k\in\supp(\phi_{uk}^*(\cI_k),\mu)=\supp(\phi_{vk}^*(\cI_k),\mu)$ and $y_k=\phi_{uk}(\overline{y}_k)=\phi_{vk}(\overline{y}_k)\in \supp(T({\cI})_k,1)$ there are affine neighborhoods $\overline{U}_k\ni \overline{y}$ and ${U}_k\ni {y}$ such that there exist  parameters $u_{1k},\ldots,u_{nk}$ and a function $v_{1k}$ on ${U}_k$  such that
\begin{enumerate}
\item $u_{1k}$ and $v_{1k}$ describe hypersurfaces of maximal contact on $U_k$ and
$\phi_{uk}^*(u_{k1})\sim \phi_{vk}^*(v_{k1})$.
\item The exceptional divisors in $E_k$ through $y_k$ are described by some parameters $u_{ik}$ and $\phi_{uk}^*(u_{ik})\sim \phi_{vk}^*(u_{ik})$.
\item $ (\phi_{u k})^*(u_{ik})-(\phi_{v k})^*(u_{ik})\in \quad {\phi}_{u 1}^*(T({\cI})_k(\overline{U}_k)={\phi}_{v k}^*(T({\cI}))_k(\overline{U}_k)$ for $i\geq 1$.
\end{enumerate}
Note then the induced marked ideals $\phi_u^*(X_i,\cH(\cI)_i,E_i,\mu)$ and $\phi_v^*(X_i,\cH(\cI)_i,E_i,\mu)$ are equal because they are controlled transforms of
$\phi_u^*(X,\cH(\cI),E,\mu)=\phi_v^*(X,\cH(\cI),E,\mu)$ defined for $(\overline{X}_i)$
(see Proposition \ref{pr: funct}).
 \qed

The above lemma can also be generalized as follows:
\begin{lemma} Let $G$ be the group of all automorphisms $\widehat{\phi}$ of 
$\widehat{X}_x$ acting trivially on the subscheme defined by $T(\cI)$ and preserving $E$, that is, $\widehat{\phi}^*(f)-f\in T(\cI)$ for any $f\in K[X_x]$ and $\phi^*(D)=D$ for any $D\in E$. Then
\begin{enumerate}
\item $\cH(\cI)$ is preserved by $G$, i.e. $\widehat{\phi}^*(\cH(\cI))_x=\cH(\cI)_x$ for any $\widehat{\phi}\in G$.
\item $G$ acts transitively on the set of tangent directions $u\in T(\cI)$ transversal to $E$.
\item Any multiple test  blow-up $(X_i)$ of $(\cI,\mu)$ is $G$-equivariant. 
\item $G$ acts trivially on the subscheme defined by $T(\cI)_i\subset X_i$.
\end{enumerate}
\end{lemma}
Since this lemma is not used in the proof and all details but a few are the same as for the proof of Lemma \ref{le: homo0} we just point out  
the differences.

\noindent{\bf Proof.}  In the proof of property (1) of Lemma \ref{le: homo0} we use the Taylor formula for $n$ unknowns 
$$\widehat{\phi}^*(f)=
f(u_1+h_1,\ldots,u_n+h_n)= f+\frac{\partial{f}}{\partial{u_1}}\cdot h_1+ \ldots +\frac{\partial{f}}{\partial{u_n}}\cdot h_n+ \frac{1}{2!}\frac{\partial^2{f}}{\partial{u_1^2}}\cdot h_1^2+
\frac{\partial^2{f}}{\partial{u_1}\partial{u_2}}\cdot h_1h_2+      \frac{1}{2!}\frac{\partial^2{f}}{\partial{u_2^2}}\cdot h_2^2+\ldots $$
In the proof of property (4) we notice that  $\widehat{\phi}$ is described by $\widehat{\phi}^*(u_i)=u_i+h_i$, where $h_i\in T(\cI)$. By a suitable linear change of coordinates we can assume that 
\begin{enumerate}
\item The center $C\subset \supp(\cI,\mu)$ of the blow-up $\sigma$ is described at $x$ by $u_1=\ldots=u_m=0$.
\item The coordinates at a point $x_1\in \sigma^{-1}(x)\cap\supp(\cI_1,\mu)$ are 
$u_{i1}=\frac{u_i}{u_m}$ for $1\leq i \leq m$ and $u_{i1}=u_i$ for $i>m$.
 \item The automorphism $\widehat{\phi}$ lifts to an automorphism $\widehat{\phi}_1$ preserving $x_1$ such that
\begin{enumerate}
\item For $i< m$, $(\widehat{\phi}_1)^*(u_{i1})=\frac{u_i+h_i}{u_m+h_m}=\frac{u_{i1}+h_i/u_m}{1+h_m/u_m}=u_{i1}+g_{i1}$ , where 
$g_{i1}\in T(\cI)_1$.
\item For $i\geq m$, $(\widehat{\phi}_1)^*(u_{i1})=u_i+h_i=u_{i1}+g_{i1}$, where $u_{i1}=u_i$ and 
$g_{i1}=h_i \in \sigma^*(T(\cI))\subset T(\cI)_1$.\qed
\end{enumerate}
\end{enumerate}

\bigskip

\subsection{Coefficient ideals and Giraud Lemma}
The idea of coefficient ideals was originated by Hironaka and then developed in papers of Villamayor and Bierstone-Milman.
The following definition modifies and generalizes the definition of Villamayor.
\begin{definition}
Let $({\cI},\mu)$ be a marked ideal of maximal order. By the
 {\it coefficient ideal}  we mean
$${\cC}({\cI},\mu)=\sum_{i=1}^\mu ({\cD}^i{\cI},\mu-i).$$
\end{definition}
\begin{remark} The coefficient ideals $\cC(\cI)$ feature two important properties.
\begin{enumerate}
\item  $\cC(\cI)$ is equivalent to  $\cI$.
\item The intersection  of the support of $(\cI,\mu)$ with any smooth subvariety $S$ is the support of the restriction of $\cC(\cI)$ to $S$:
$$\supp(\cI)\cap S=\supp(\cC(\cI)_{|S}).$$ 
\noindent Moreover this condition is persistent under relevant multiple test  blow-ups.
\end{enumerate}
These properties allow one to control and modify the part of support of $(\cI,\mu)$ contained in $S$ by applying multiple test  blow-ups of  $\cC(\cI)_{|S}$.
\end{remark}

\begin{lemma}\label{le: coeff}
 ${\cC}({\cI},\mu)\simeq({\cI},\mu)$.

\end{lemma}
\noindent {\bf Proof.} By Lemma \ref{le: operations} multiple test  blow-ups of $\cC(\cI,\mu)$ are simultaneous multiple test  blow-ups of $\cD^i(\cI,\mu)$ for $0\leq i\leq \mu-1$. By Lemma \ref{le: inclusions2} multiple test  blow-ups of $(\cI,\mu)$ define the multiple test  blow-up of all $\cD^i(\cI,\mu)$.
Thus multiple test  blow-ups of $(\cI,\mu)$ and $\cC(\cI,\mu)$ are the same and $\supp({\cC}({\cI},\mu))_k=\bigcap \supp({\cD}^i{\cI},\mu-i)_k=\supp({\cI}_k,\mu).$
\qed

\begin{lemma} \label{le: cf} Let $(X,{\cI},E,\mu)$ be a marked ideal of maximal order  whose support $\supp({\cI},\mu)$ does not contain a smooth subvariety $S$ of $X$. Assume that $S$ has only simple normal crossings with $E$. Then 
$$\supp({\cI},\mu)\cap S\subseteq \supp(({\cI},\mu)_{|S}).$$ 
Let $\sigma: X'\to X$ be a blow-up with center $C\subset \supp({\cI},\mu)\cap S$. Denote by $S'\subset X'$ the strict transform of $S\subset X$. Then
$$\sigma^{\rm c}((\cI,\mu)_{|S})=(\sigma^{\rm c}(\cI,\mu))_{|S'}.$$
Moreover for any  multiple test blow-up $(X_i)$   with all centers $C_i$   contained in the strict transforms $S_i\subset X_i$ of $S$,
the restrictions  $\sigma_{i|S_i}: S_i\to S_{i-1}$ of the morphisms $\sigma_i: X_i\to X_{i-1}$ to $S_i$ define
  a multiple test  blow-up $(S_i)$
of $({\cI},\mu)_{|S}$ such that
$$[(\cI,\mu)_{|S}]_i=(\cI_i,\mu)_{|S_i}.$$

\end{lemma}

\noindent{\bf Proof}. The first inclusion holds since the order of an ideal does not drop but may rise after restriction to a subvariety. Let $x_1,\ldots,x_k$ describethe subvariety $S$ of $X$ at a point $p\in C$. Let $p'\in S'$ map to $p$. We can
find coordinates $x_1,\ldots,x_k,y_1,\ldots,y_{n-k}$ such that 
the center of the blow-up is described by $x_1,\ldots,x_k,y_1,\ldots,y_{m}$ and the 
coordinates at $p'$ are given by $$x'_1=x_1/y_m, \ldots, x'_k=x_k/y_m, 
y'_1= y_1/y_m,\ldots, y'_{m}=y_m,y'_{m+1}=y_{m+1},\ldots, y'_n=y_n$$   
\noindent where the strict transform $S'\subset X'$ of $S$ is described by $x'_1,\ldots,x'_k$.
Then we can write a  function $f\in \cI(U)$ as $f=\sum c_{\alpha f}(y) {x}^\alpha$,  where $c_{\alpha f}(y)$ are formal power series in $y_i$. The controlled transform $f'=\sigma^{\rm c}(f)=y^{-\mu}(f\circ\sigma)$ can be written as
$$f'=\sum c'_{\alpha f}(y) {x'}^\alpha, $$
\noindent where $c'_{\alpha f}=y_m^{-\mu+|\alpha|} \sigma^*(c_{\alpha f})$.
But then $f_{|S}=(c_{0 f})_{|S}$ and $$\sigma^{\rm c}(f)_{|S'}=(c'_{0 f})_{|S'}=y_m^{-\mu}\sigma^*(c_{0 f})_{|S'}=y_m^{-\mu}\sigma^*((c_{0 f})_{|S})=y_m^{-\mu}\sigma^*(f_{|S})=\sigma^{\rm c}(f_{|S}).$$
The last  part of the theorem follows by induction:
\begin{eqnarray}
&\supp({\cI}_{i+1},\mu)\cap S_{i+1}=\supp(\sigma_{i+1}^{\rm c}({\cI}_{i},\mu))\cap S_{i+1}\subseteq \supp(\sigma_{i+1|S_i})^{\rm c}(({\cI}_{i},\mu)_{|S_i}) \subseteq \nonumber \\ &\supp((\sigma_{i+1|S_i})^{\rm c}({\cI},\mu)_{|S})_i=\supp(({\cI},\mu)_{|S})_{i+1},\nonumber
\end{eqnarray}

$$[(\cI,\mu)_{|S}]_{i+1}=\sigma_{i+1|S_i}^{\rm c}[(\cI,\mu)_{|S}]_i=\sigma_{i+1|S_i}^{\rm c}((\cI_i,\mu)_{|S_i})=(\sigma_{i+1}^{\rm c}(\cI_i,\mu))_{|S_{i+1}}=(\cI_i,\mu)_{|S_{i+1}}.$$
\qed

\begin{lemma} \label{le: S} Let $(X,{\cI},E,\mu)$ be a marked ideal of maximal order  whose support $\supp({\cI},\mu)$ does not contain a smooth subvariety $S$ of $X$. Assume that $S$ has only simple normal crossings with $E$. Then 
$$\supp({\cI},\mu)\cap S= \supp({\cC}({\cI},\mu)_{|S}).$$ 
Moreover  let  $(X_i)$ be  a multiple test blow-up with  centers $C_i$   contained in the strict transforms $S_i\subset X_i$ of $S$. Then

\begin{enumerate}

\item The restrictions  $\sigma_{i|S_i}: S_i\to S_{i-1}$ of the morphisms $\sigma_i: X_i\to X_{i-1}$ define
  a multiple test  blow-up $(S_i)$
of ${\cC}({\cI},\mu)_{|S}$.

\item $\supp({\cI}_i,\mu)\cap S_i= \supp[{\cC}({\cI},\mu)_{|S}]_i.$
\item Every multiple test  blow-up $(S_i)$  of ${\cC}({\cI},\mu)_{|S}$  defines  a multiple test  blow-up $(X_i)$ of $({\cI},\mu)$ with  centers $C_i$ contained in the strict transforms $S_i\subset X_i$ of $S\subset X$.

\end{enumerate}

\end{lemma}

\noindent {\bf Proof.}
By Lemmas \ref{le: coeff} and \ref{le: cf}, $\supp({\cI},\mu)\cap S= \supp ({\cC}({\cI},\mu))\cap S\subseteq\supp({\cC}({\cI},\mu)_{|S})
$.

 Let $x_1,\ldots,x_k,y_1,\ldots, y_{n-k}$ be local parameters at $x$ such that $\{x_1=0,\ldots,x_k=0\}$ describes $S$. Then any function $f\in {\cI}$ can be written as 
$$f=\sum c_{\alpha f}(y) x^\alpha,$$\noindent where $c_{\alpha f}(y)$ are formal power series in $y_i$.

Now $x\in \supp({\cI},\mu)\cap S$ iff $\ord_x(c_\alpha)\geq \mu-|\alpha|$ for all $f\in {\cI}$ and $|\alpha|\leq \mu$. Note that $$c_{\alpha f|S}=\bigg(\frac{1}{\alpha!}\frac{\partial^{|\alpha|}(f)}{\partial x^\alpha}\bigg)_{|S}\in {\cD}^{|\alpha|}({\cI})_{|S}$$ and consequently  $\supp({\cI},\mu)\cap S=\bigcap_{f\in {\cI}, |\alpha|\leq \mu}\supp({c_{\alpha f|S}}, \mu-|\alpha|)\supset  \supp({\cC}({\cI},\mu)_{|S})$.

Assume that all multiple test  blow-ups of $({\cI},\mu)$ of length $k$ with centers $C_i\subset S_i$ are defined by multiple test  blow-ups of ${\cC}({\cI},\mu)_{|S}$ and moreover for $i\leq k$,
$$\supp({\cI}_i,\mu)\cap S_i=\supp[{\cC}({\cI},\mu)_{|S}]_i.$$
For any $f\in \cI$ define $f=f_0\in\cI$ and $f_{i+1}=\sigma^{\rm c}_i(f_i)=y_i^{-\mu}\sigma^*(f_i)\in \cI_{i+1}$.
Assume moreover that for any $f\in {\cI}$, $$f_{k}=\sum c_{\alpha fk}(y) x^\alpha ,$$\noindent where $c_{\alpha f{k}|S_{k}}\in {(\sigma^{k}_{|S_{k}})^{\rm c}}({\cD}^{\mu-|\alpha|}({\cI})_{|S})$. 
Consider the effect of the blow-up of $C_k$ at a point $x'$ in the strict transform of $S_{k+1}\subset X_{k+1}$.
By Lemmas \ref{le: coeff} and \ref{le: cf},
$$\supp({\cI}_{k+1},\mu)\cap S_{k+1}= \supp [{\cC}({\cI},\mu)]_{k+1}\cap S_{k+1}\subseteq\supp[{\cC}({\cI},\mu)]_{{k+1}|S_{k+1}}=\supp[{\cC}({\cI},\mu)_{|S}]_{k+1}
$$

Let $x_1,\ldots,x_k$ describe the subvariety $S_k$ of $X_k$. We can
find coordinates $x_1,\ldots,x_k,y_1,\ldots,y_{n-k}$, by taking if necessary linear combinations of $y_1,\ldots,y_{n-k}$, such that 
the center of the blow-up is described by $x_1,\ldots,x_k,y_1,\ldots,y_{m}$ and the 
coordinates at $x'$ are given by $$x'_1=x_1/y_m, \ldots, x'_k=x_k/y_m, 
y'_1= y_1/y_m,\ldots, y'_{m}=y_m,y'_{m+1}=y_{m+1},\ldots, y'_n=y_n.$$   
Note that replacing $y_1,\ldots,y_{n-k}$ with their linear combinations does not modify the form $f_{k}=\sum c_{\alpha fk}(y) x^\alpha .$
Then the function $f_{k+1}=\sigma^{\rm c}(f_k)$ can be written as
$$f_{k+1}=\sum c_{\alpha f k+1}(y) {x'}^\alpha, $$
\noindent where $c_{\alpha f{k+1}}=y_m^{-\mu+|\alpha|} \sigma_{k+1}^*(c_{\alpha f k})$. Thus $$c_{\alpha f{k+1}|S_{k+1}}={(\sigma_{k+1|S_{k+1}})^{\rm c}}(c_{\alpha f{k}|S_{k}})\in {(\sigma^{k+1}_{|S_{k+1}})^{\rm c}}({\cD}^{\mu-|\alpha|}({\cI})_{|S})={(\sigma^{k+1})^{\rm c}}({\cD}^{\mu-|\alpha|}({\cI}))_{|S_{k+1}}$$ and consequently 
 $$\supp({\cI}_{k+1},\mu)\cap S_{k+1}=\bigcap_{f\in {\cI}, |\alpha|\leq \mu}
 \supp({c_{{\alpha f{k+1}} |S_{k+1}}}, \mu-|\alpha|)\supseteq  \supp[{\cC}({\cI},\mu)_{|S}]_{k+1}=\supp({\cC}({\cI},\mu)_{k+1})_{|S_{k+1}}.$$ 
\qed 

A direct consequence of the above lemma is the following result:
\begin{lemma} \label{le: S2} Let $(X,{\cI},E,\mu)$ be a marked ideal of maximal order  whose support $\supp({\cI},\mu)$ does not contain a smooth subvariety $S$ of $X$. Assume that $S$ has only simple normal crossings with $E$. Let $(X_i)$ be its multiple test 
 blow-up such that all centers $C_i$  are either contained in the strict transforms $S_i\subset X_i$ of $S$ or are disjoint from them.
Then the restrictions  $\sigma_{i|S_i}: S_i\to S_{i-1}$ of the morphisms $\sigma_i: X_i\to X_{i-1}$ define
  a multiple test  blow-up $(S_i)$ of ${\cC}({\cI},\mu)_{|S}$ and 
\begin{equation}\nonumber
\supp({\cI}_i,\mu)\cap S_i= \supp[{\cC}({\cI},\mu)_{|S}]_i.
\end{equation}
\end{lemma}

As a simple consequence of the Lemma \ref{le: S} we formulate the following   refinement of the Giraud Lemma.

\begin{lemma}\label{le: codim2} Let $(X,{\cI},\emptyset,\mu)$ be a marked ideal of maximal order  whose support $\supp({\cI},\mu)$ has codimension at least $2$ at some point $x$.
Let $U\ni x$ be an open subset for which there is
a tangent direction $u\in T({\cI})$ and such that  $\supp({\cI},\mu)\cap U$ 
is of codimension 2. 
Let $V(u)$ be the regular subscheme of $U$ defined by $u$. Then for any multiple test blow-up  $X_i$ of
$X$,
\begin{enumerate}
\item $\supp({\cI}_i,\mu)$ is contained in the strict transform $V(u)_i$ of $V(u)$ as a proper subset.
\item The sequence  $(V(u)_i)$ is  a multiple test  blow-up  of ${\cC}({\cI},\mu)_{|V(u)}$.

\item $\supp({\cI}_i,\mu)\cap V(u)_i= \supp[{\cC}({\cI},\mu)_{|V(u)}]_i.$
\item Every multiple test blow-up $(V(u)_i)$  of ${\cC}({\cI},\mu)_{|V(u)}$ defines a multiple test  blow-up  $(X_i)$ of $({\cI},\mu)$.
\end{enumerate}\qed
\end{lemma}

\begin{lemma} \label{le: etale3} Let $\phi: X'\to X$ be a smooth morphism of smooth varieties 
and let $(X,{\cI},\emptyset,\mu)$ be a marked ideal. Then 

$$\phi^*({\cC}({\cI}))={\cC}(\phi^*({\cI})).$$
\end{lemma}
\noindent{\bf Proof.} A direct consequence of Lemma \ref{le: etale}. \qed
\bigskip

\section{Resolution algorithm}
The presentation of the following Hironaka resolution algorithm  builds upon 
Villamayor's and Bierstone-Milman's proofs.

Let $\Sub(E_i)$ denote the set of all divisors of $E_i$.  For any  subset in $\Sub(E_i)$ write a sequence consisting of all elements of the subset ordered lexicographically followed by the infinite sequence of zeros $({D}_1,{D}_2,\ldots,0,\ldots)$. We shall assume that $0\leq {D}$ for any ${D}$.
Then for any two subsets $A_1=\{{D}^1_i\}_{i\in I}$ and $A_2=\{{D}^2_j\}_{j\in {\cJ}}$ for which $\{{D}^1_i\}_{i\in I}\cap \{{D}^2_j\}_{j\in {\cJ}}\neq\emptyset$  we write
$$A_1\leq A_2$$ \noindent if
for the corresponding sequences $(D^1_1,D^1_2,\ldots,0,\ldots)\leq (D^2_1,D^2_2,\ldots,0,\ldots)$. 

Let ${\bf Q}_{\geq 0}$ denote the set of nonnegative rational numbers and let
$$\overline{\bf Q}_{\geq 0}:= {\bf Q}_{\geq 0}\cup\{\infty\}.$$
 Denote by $\overline{\bf Q}^\infty_{\geq 0}$ the set of all infinite sequences in $\overline{\bf Q}_{\geq 0}$ with a finite number of nonzero elements. We equip $\overline{\bf Q}^\infty_{\geq 0}$ with the lexicographical order.

\begin{proposition}\label{pr: resol}  For any  marked ideal $(X,{\cI},E,\mu)$ such that $\cI\neq 0$ there is an associated  resolution $(X_i)_{0\leq i\leq m_X}$, called \underline{canonical},
 satisfying the following conditions:

\begin{enumerate}
\item There are  upper semicontinuous invariants $\inv$, $\nu$ and $\rho$ defined on  $\supp(X_i,{\cI}_i,E_i,\mu)$ with values in ${\bf Q}_{\geq 0}\times\overline{\bf Q}^\infty_{\geq 0}$, $\overline{\bf Q}_{\geq 0}$ and  $\Sub({E_i})$ respectively. 
\item The centers $C_i$ of blow-ups are  regular and defined by the set where 
$(\inv,\rho)$ attains its maximum. They are components of  the maximal locus of $\inv$. 
\item
\begin{enumerate} 
\item For any $x\in \supp(X_{i+1},{\cI}_{i+1},E_{i+1},\mu)$ and $\sigma(x)\in
C_i$, either $\inv(x)<\inv(\sigma(x))$ or \\$\inv_{i+1}(x)=\inv(\sigma(x))$ and  $\mu(x)<\mu(\sigma(x))$. 

\item For any $x\in \supp(X_{i+1},{\cI}_{i+1},E_{i+1},\mu)$ and $\sigma(x)\not\in
C_i$,  $\inv(x)=\inv(\sigma(x))$, $\rho(x)=\rho(\sigma(x))$ and  $\mu(x)=\mu(\sigma(x))$. 
\end{enumerate}

\item For any  \'etale morphism $\phi: X'\to X$ the induced sequence 
 $(X'_i)=\phi^*(X_i)$ is an extension   of the canonical resolution of $X'$ such that for the induced marked ideals $(X'_i,\cI'_i,E'_i,\mu)$ and  $x'\in\supp(X'_i,\cI'_i,E'_i,\mu)$, we have $\inv(\phi_i(x'))=\inv(x')$, $\nu(\phi_i(x'))=\nu(x')$ and $\rho(\phi_i(x'))=\rho(x')\in \Sub(E'_i)\subset \Sub(E_i)$.

\end{enumerate}

\end{proposition}

\noindent {\bf Proof.} Induction  on the dimension of $X$. If $X$ is $0$-dimensional, $\cI\neq 0$ and $\mu>0$ then $\supp(X, {\cI},\mu)=\emptyset$ and all resolutions are trivial.

The invariants $\inv$, $\nu$ and $\rho$ will be defined successively in the course of the resolution algorithm using the inductive assumptions and property (2).

{\bf Step 1}.  {\bf Resolving a marked ideal $(X,{\cJ},E,\mu)$ of maximal
order.}

The process of resolving  the marked ideals of maximal order is controlled by an auxilary invariant $\overline{\inv}$ defined in Step 1. The invariant $\inv$ will then be defined for any marked ideals in  Step 2.

Before we start our  resolution algorithm for the marked ideal $({\cJ},\mu)$ of maximal order we shall replace it with the equivalent  homogenized ideal ${\cC}({\cH}({\cJ},\mu))$. 
Resolving the ideal ${\cC}({\cH}({\cJ},\mu))$ defines a resolution of $({\cJ},\mu)$ at this step.
To simplify notation we shall denote ${\cC}({\cH}({\cJ},\mu))$ by $(\overline{\cJ},\overline{\mu})$.
 
 {\bf Step 1a.}  {\bf Reduction to the nonboundary case.}
 For any multiple test  blow-up $(X_i)$ of $(X,{\overline{\cJ}},E,\overline{\mu})$ we shall identify (for simplicity) strict transforms of $E$ on $X_i$ with $E$. 

For any $x\in X_i$, let $s(x)$ denote the number of divisors in $E$ through $x$ and 
set $$ s_i=\max\{ s(x) \mid x \in \supp({\overline{\cJ}}_i)\}.$$

Let $s=s_0$. By  assumption the intersections of any $s> s_0$ components of the exceptional divisors are disjoint from 
$\supp({\overline{\cJ}},\overline{\mu})$.
Each intersection of divisors in $E$ is locally defined by intersection of some irreducible components of these divisors.
Find all   intersections $H^s_\alpha, \alpha\in A$, of $s$ irreducible components of divisors $E$ such that $\supp({\overline{\cJ}},\overline{\mu})\cap H^s_\alpha\neq\emptyset$. By the maximality of $s$, the supports $\supp({\overline{\cJ}}_{|H^s_\alpha})\subset H^s_\alpha$ are disjoint from $H^s_{\alpha'}$, where $\alpha'\neq\alpha$.

\quad {\bf Step 1aa.} {\bf Eliminating the components $H^s_\alpha$ contained in $\supp({\overline{\cJ}},\overline{\mu})$.}

Let $H^s_\alpha \subset \supp({\overline{\cJ}},\overline{\mu})$. If $s\geq 2$ then by blowing up  $C=H^s_a$ we separate divisors contributing to $H^s_a$, thus creating new points all with $s(x)<s$. 
If $s=1$  then by Lemma \ref{le: codim1}, $H^s_\alpha \subset \supp({\overline{\cJ}},\overline{\mu})$ is a codimension one component and by blowing up $H^s_\alpha$ we create all new points off $\supp({\overline{\cJ}},\overline{\mu})$.

For all $x\in H^s_a\subset \supp({\overline{\cJ}},\overline{\mu})$ set 
$$\overline{\inv}(x) =(s(x),\infty,0,\ldots,0,\ldots),\quad \nu(x)=0,\quad  \rho(x)=\emptyset.$$ \noindent 
This definition, as we see below, is devised so as to ensure that all $H^s_\alpha\subset \supp({\overline{\cJ}},\overline{\mu})$ will be blown up first and we reduce the situation to the case where no  $H^s_\alpha$ is contained in $\supp({\overline{\cJ}},\overline{\mu})$.

\quad {\bf Step 1ab.} {\bf Moving $\supp({\overline{\cJ}},\overline{\mu})$ and $H^s_\alpha$ apart .}

After the blow-ups in Step 1aa we arrive at $X_p$  for which no $H^s_\alpha$ is contained in $\supp({\overline{\cJ}}_p,\overline{\mu})$, where $p=0$ if there were no such components and $p=1$ if there were some. 

Construct the canonical resolutions of ${\overline{\cJ}_p}_{|H^s_\alpha}$. By Lemma \ref{le: S} each such resolution 
defines a multiple test  blow-up of $({\overline{\cJ}}_p,\overline{\mu})$ (and of $({\overline{\cJ}},\overline{\mu})$) . Since the supports $\supp({\overline{\cJ}}_{|H^s_\alpha})\subset H^s_\alpha$ are disjoint from $H^s_{\alpha'}$, where $\alpha'\neq\alpha$,  these resolutions
glue to a unique multiple test  blow-up $(X_i)_{i\leq j_1}$ of $({\overline{\cJ}},\overline{\mu})$ such that  $s_{j_1}<s$. 
To control the glueing procedure and ensure its uniqueness  we  define for all $x\in \supp({\overline{\cJ}},\overline{\mu})\cap H^s_\alpha$ the invariant $$ \overline{\inv}(x)=(s(x),\inv_ {{\overline{\cJ}}_{|H^s_\alpha}}(x)), \quad\nu_{\overline{\cJ}}=\nu_{{\overline{\cJ}}_{|H^s_\alpha}},\quad \rho_{\overline{\cJ}}=\rho_{{\overline{\cJ}}_{|H^s_\alpha}}. $$ \noindent The blow-ups will be performed at the centers $C\subset \supp({\overline{\cJ}},\overline{\mu})\cap H^s_\alpha$ for which the invariant $(\overline{\inv},\rho)$ attains its maximum. 
Note that by the maximality condition for any $H^s_\alpha$ the irreducible components of the centers are contained in $H^s_\alpha$ or are disjoint from them. Therefore by Lemma \ref{le: S2}, $$\supp(\overline{\cJ}_i,\overline{\mu})_{|H^s_\alpha}=\supp(\overline{\cJ}_i,\overline{\mu})\cap{H^s_\alpha}.$$
By applying this multiple test blow-up we create a marked ideal $({\overline{\cJ}}_{j_1},\overline{\mu})$  with support disjoint from all  $H^s_\alpha$.
Summarizing the above we construct  a multiple test blow-up $(X_i)_{0\leq i\leq  j_1}$ subject to the conditions:
\begin{enumerate}
\item  $(H^s_{\alpha i})_{0\leq i\leq j_1}$ is an extension  of the canonical resolution  of ${\overline{\cJ}}_{|H^s_{\alpha}}$.

\item There are invariants $\overline{\inv}$, $\overline{\mu}$ and $\rho$ defined for 
${0\leq i< j_1}$ and all
$x\in \supp({\overline{\cJ}_i},\overline{\mu})\cap H^s_{\alpha i}$ such that  $$ \overline{\inv}(x)=(s(x),\inv_ {{\overline{\cJ}}_{i|H^s_\alpha}}(x)), \quad\nu_{\overline{\cJ}_i}=\nu_{{\overline{\cJ}}_{i|H^s_\alpha}},\quad \rho_{\overline{\cJ_i}}=\rho_{{\overline{\cJ}}_{i|H^s_\alpha}}. $$ 
\item The blow-ups of $X_i$ are performed at the centers where the invariant $(\overline{\inv},\rho)$ attains its maximum.
\item $\supp({\overline{\cJ}_{j_1}},\overline{\mu})\cap H^s_{\alpha j_1}=\emptyset$.
 \end{enumerate}

\quad {\bf Conclusion of the algorithm in Step 1a.} 
After performing the blow-ups in Steps 1aa and 1ab for the marked ideal $(\overline{\cJ},\overline{\mu})$ we arrive at a marked ideal $({\overline{\cJ}}_{j_1},\overline{\mu})$ with $s_{j_1}<s_0$. Now we put $s=s_{j_1}$ and repeat the procedure of Steps 1aa and 1ab   for
 $({\overline{\cJ}}_{j_1},\overline{\mu})$. Note that  any $H^s_{\alpha j_1}$ on $X_{j_1}$ is the strict transform of some intersection $H^{s_{j_1}}_\alpha$ of $s=s_{j_1}$ divisors in $E$ on $X$. Moreover by the maximality condition for all $s_i$, where $i\leq {j_1}$ and $\alpha\neq\alpha'$, the set $\supp(\overline{\cJ}_i,\overline{\mu})\cap H^{s_i}_{\alpha' i}$ is either disjoint from $H^{s_{j_1}}_{\alpha i}$ or  contained in it. Thus for $0\leq i\leq {j_1}$, all centers $C_i$ have components either contained in  $H^{s_{j_1}}_{\alpha i}=H^s_{\alpha i}$ or  disjoint from them and by Lemma \ref{le: S2}, \begin{equation} \supp(\overline{\cJ}_{i},\overline{\mu})_{|H^s_{\alpha i}}=\supp(\overline{\cJ}_{i},\overline{\mu})\cap{H^s_{\alpha i}}.\nonumber
 \end{equation}
Moreover if we repeat the procedure in Steps 1aa and 1ab 
the above property will still be satisfied until either $(\overline{\cJ}_{i},\overline{\mu})_{|H^s_\alpha}$ are resolved as in Step 1ab or  $H^s_\alpha$ disappear as in Step 1aa. 

We continue the above process till $s_{j_k}=s_r=0$. 
Then $(X_j)_{0\leq j\leq r}$  is a multiple test  blow-up of $(X,{\overline{\cJ}},E,\overline{\mu})$
such that $\supp({\overline{\cJ}}_r,\overline{\mu})$ does not intersect any divisor in $E$.
Therefore $(X_j)_{0\leq j\leq r}$ and further longer multiple test  blow-ups $(X_j)_{0\leq j\leq r_0}$ for any $r\leq r_0$   can be considered as multiple test  blow-ups of $(X,{\overline{\cJ}},\emptyset,\overline{\mu})$ since starting from $X_r$ the strict transforms of $E$ play no further role in the resolution process since they do not intersect  $\supp({\overline{\cJ}}_j,\overline{\mu})$ for $j\geq r$.

 Note that in Step 1a all points $x\in\supp(\overline{\cJ}_i,\overline{\mu})$ 
for which $s(x)>0$ were assigned their invariants $\overline{\inv}$, $\nu$ and $\rho$.
(They are assigned the invariants at the moment they are getting blown-up. The invariants remain unchanged when the points are transformed isomorphically.)
The invariants are upper semicontinuous by the semicontinuity of the function $s(x)$ and the inductive assumption.

\bigskip
{\bf Step 1b}.   {\bf Nonboundary case}

Let $(X_j)_{0\leq j\leq r}$ be the multiple test  blow-up of $(X,{\overline{\cJ}},\emptyset,\overline{\mu})$ defined in Step 1a.

\quad {\bf Step 1ba}. {\bf Eliminating the codimension one components of $\supp({\overline{\cJ}}_r,\overline{\mu})$.}

If $\supp({\overline{\cJ}}_r,\overline{\mu})$ is of codimension $1$ then by Lemma \ref{le: codim1} all its codimension $1$ 
components are smooth and disjoint from the other components of $\supp({\overline{\cJ}}_r,\overline{\mu})$. These components are strict transforms of the codimension $1$ components of $\supp({\overline{\cJ}},\overline{\mu})$. Moreover the irreducible components of the centers of blow-ups were either contained in the strict transforms or disjoint from them. Therefore $E_r$ will be transversal to all the codimension $1$ components. 
Let $\codim(1)(\supp({\overline{\cJ}}_i,\overline{\mu}))$ be the union of all components of $\supp({\overline{\cJ}}_i,\overline{\mu}))$ of codimension $1$. We define the invariants for $x\in \codim(1)(\supp({\overline{\cJ}}_r,\overline{\mu})$ to be $$\overline{\inv}(x) =(0,\infty,0,\ldots,0,\ldots),\quad \nu(x)=0,\quad  \rho(x)=\emptyset.$$ \noindent 
This definition, as we see below, is devised so as to ensure that all codimension $1$ components will be blown up first. 

By Lemma \ref{le: codim1} blowing up the components reduces the situation to the case when  $\supp({\overline{\cJ}},\overline{\mu})$ is of codimension $\geq 2$.

\quad {\bf Step 1bb}. {\bf Eliminating the codimension $\geq 2$ components of $\supp({\overline{\cJ}}_r,\overline{\mu})$.}

After  Step 1ba we arrive  at a marked ideal $\supp({\overline{\cJ}}_p,\overline{\mu})$,
where $p=r$ if there were no codimension one components and $p=r+1$ if there were 
some and we blew them up.

For any $x\in \supp({\overline{\cJ}},\overline{\mu})\setminus \codim(1)(\supp({\overline{\cJ}},\overline{\mu})\subset X $ find a tangent direction $u\in {\cD}^{\overline{\mu}-1}({\overline{\cJ}})$  on some neighborhood $U_u$ of $x$. 
Then $V(u)\subset U_u$ is a hypersurface of maximal contact.
By the quasicompactness of $X$ we can assume that the covering defined by $U_u$ is finite.
Let $U_{ui}\subset X_i$ be the inverse image of $U_u$ and let $V(u)_i\subset U_u$ denote the strict transform of $V(u)$. By Lemma \ref{le: codim2},
$(V(u)_i)_{0\leq i\leq p}$ is a  multiple test  blow-up of $(V(u), {\overline{\cJ}}_{|V(u)},\emptyset,\overline{\mu})$. In particular the induced marked
 ideal for $i=p$ is equal to
$${\overline{\cJ}}_{p|V(u)_{p}}=(V(u)_p,{\overline{\cJ}}_{p|V(u)_{p}}, (E_p\setminus E)_{|V(u)_{p}},\overline{\mu}).$$
Construct the canonical resolution of $(V(u)_i)_{p\leq i\leq m_u}$ of the marked ideal ${\overline{\cJ}}_{p|V(u)_{p}}$.
\noindent Then the sequence  $(V(u)_i)_{0\leq i \leq m_u}$ is a resolution of  $(V(u), {\overline{\cJ}}_{|V(u)},\emptyset,\overline{\mu})$ which  defines, 
by Lemma \ref{le: codim2}, a resolution $(U_{ui})_{0\leq i\leq m_u}$
of $(U_u,{\overline{\cJ}}_{|U_u},\emptyset, \overline{\mu})$. Moreover both resolutions are related by the property $$\supp(\overline{\cJ}_{i|U_{ui}})=\supp(\overline{\cJ}_{i|V(u)_i}).$$

We shall construct the resolution of $(X,{\overline{\cJ}},\emptyset,\overline{\mu})$  by patching together extensions of the local resolutions $(U_{ui})_{0\leq i\leq m_u}$.

For  $x\in \supp(\overline{\cJ}_p,\overline{\mu})\cap U_{up}$ define the invariants $$\overline{\inv}(x):=(0,\inv_{{\overline{\cJ}}_{p|V(u)_p}}(x)),\quad \nu:=\nu_{{\overline{\cJ}}_{p|V(u)_p}}(x),\quad\rho(x):=
\rho_{{\overline{\cJ}}_{p|V(u)_p}}(x).$$  
We need to show that these invariants do not depend on the choice of $u$.

Let $x\in \supp(\overline{\cJ}_p,\overline{\mu})\cap U_{up}\cap U_{vp}$. 
 By Glueing Lemma \ref{le: homo} for any two different tangent directions $u$ and $v$ we find  \'etale neighborhoods $\phi_u,\phi_v: U^{uv}\to U:=U_u\cap U_v$ and their liftings 
$\phi_{pu},\phi_{pv}: U_p^{uv}\to U_p:=U_{up}\cap U_{vp}$ such that 
\begin{enumerate}
\item $X_p^{uv}:=(\phi_{pu})^{-1}(V(u)_p)=(\phi_{pv})^{-1}(V(v)_p)$.
\item  $(U^{uv}_p , \overline{\cJ}^{uv}_p, E^{uv}_p,\overline{\mu}):=(\phi_{pu})^{*}(U_p,\overline{\cJ}_p,E_p,\overline{\mu})=(\phi_{pv})^{*}(U_p,\overline{\cJ}_p,E_p,\overline{\mu})$.   
\item There exists $y\in \supp(U^{uv}_p , \overline{\cJ}^{uv}_p, E^{uv}_p,\overline{\mu})$
such that
$\phi_{pu}(\overline{x})=\phi_{pv}(\overline{x})$.
\end{enumerate}
 By the functoriality  of the invariants we have
\begin{equation}\inv_{{\overline{\cJ}}_{p|V(u)_p}}(x)\quad = \quad \inv_{ {\overline{\cJ}}^{uv}_{p\,  |X^{uv}_p}}(\overline{x})\quad = \quad \inv_{{\overline{\cJ}}_{p|V(v)_p}}(x).\nonumber\end{equation}
Analogously $\nu_{{\overline{\cJ}}_{p|V(u)_p}}(x)=\nu_{{\overline{\cJ}}_{p|V(v)_p}}(x)$ \quad  and \quad $\rho_{{\overline{\cJ}}_{p|V(u)_p}}(x)=\rho_{{\overline{\cJ}}_{p|V(v)_p},}(x)$.
Thus the invariants $\overline{\inv}$, $\nu$ and $\rho$ do not depend on the choice of the tangent direction.

Define the center $C_p$ of the  blow-up on $X_p$  to be the maximal locus of the invariant $(\overline{\inv},\rho)$. 
Note that for any tangent direction $u$, either $C_p\cap U_{up}$ defines the first blow-of the canonical resolution of \\$(V(u)_{p},{\overline{\cJ}}_{p|V(u)_{p}},E_{p|V(u)_{p}},\overline{\mu})$ or $C_p\cap U_{up}=\emptyset$ and the blow-up of $C_p$ does not change $V(u)_p\subset U_{up}$.

Blowing up $C_p$ defines $X_{p+1}$ and we are in a position to construct the invariants on $X_{p+1}$ and define the center of the blow-up $C_{p+1}\subset X_{p+1}$ as before.

By repeating the same reasoning  for   $j=p+1, \ldots, m$ we construct the resolution $(X_i)_{p\leq i\leq m}$ of \\ $(X_p,{\overline{\cJ}}_p,E_p\setminus E,\overline{\mu})$ satisfying the following properties.

\begin{enumerate}
\item For any $u$,  the restriction of $(X_i)_{p\leq i\leq m}$  to $(V(u)_i)_{ p\leq i\leq m}$ is an extension  of the canonical resolution of  $(V(u)_{p},{\overline{\cJ}}_{p|V(u)_{p}},E_{p|V(u)_{p}},\overline{\mu})$.

\item There are invariants $\overline{\inv}$, $\overline{\mu}$ and $\rho$ defined for  all points
$x\in\supp({\overline{\cJ}}_i,\overline{\mu})$, ${ p\leq i\leq m}$,  such that $$\overline{\inv}(x):=(0,\inv_{{\overline{\cJ}}_{i|V(u)_i}}),\quad \nu(x):=\nu_{{\overline{\cJ}}_{i|V(u)_i}}(x),\quad\rho(x):=
\rho_{{\overline{\cJ}}_{i|V(u)_i}}(x).$$ 

\item The blow-ups of $X_i$ are performed at the centers where the invariant $(\overline{\inv},\rho)$ attains its maximum.

\item $\supp(\overline{\cJ}_m,\overline{\mu})=\emptyset.$
 \end{enumerate}

The resolution $(X_i)_{p\leq i\leq m}$ of $(X_p,{\overline{\cJ}}_p,E_p\setminus E,\overline{\mu})$  defines the resolution $(X_i)_{0\leq i\leq m}$ of $(X,{\overline{\cJ}},\emptyset,\overline{\mu})$  and of $(X,{\overline{\cJ}},E,\overline{\mu})$.

In  Step 1b all points $x\in \supp(\overline{\cJ}_i,\overline{\mu})$ with $s(x)=0$ were assigned the invariants $\overline{\inv}$, $\nu$ and $\rho$. They are upper semicontinuous 
by the inductive assumption.

\bigskip
{\bf Commutativity of the resolution procedure in Step 1 with \'etale morphisms.}

Let $\phi: X'\to X$ be an \'etale morphism. 
In Step 1a we find a sequence $i_0:=0<i_1<\ldots<i_k=r\leq m$ such that $s_{i_0}>s_{i_1}>\ldots>s_{i_k}$ and for  ${i_l}\leq i <i_{l+1}$, we have $s_i=s_{i_l}$.
Moreover the resolution process for $(X_i)_{i_l\leq i\leq i_{l+1}}$ is reduced to resolving  $\overline{\cJ}_{i_{l}|H^s_{\alpha i_l}}$. In Step 1b we reduce the resolution process for
$(X_i)_{i_k\leq i\leq m}$ to resolving $\overline{\cJ}_{i_k|V(u)_{i_k}}$.

 Let $s'_{j_0}>s'_{j_1}>\ldots>s'_{j_{k'}}$ be the  corresponding sequence defined for the canonical resolution $(X'_i)_{0\leq j\leq m'}$ of $$(X',\overline{\cJ}',E',\overline{\mu}):=\phi^*(X,\overline{\cJ},E,\overline{\mu}).$$ Let $\phi^*(X_i)_{0\leq i\leq m}$ denote the resolution  of $(\overline{\cJ'},\overline{\mu})$ 
induced by $(X_i)_{0\leq i\leq m}$. In particular $X'_0=\phi^*(X_0)$.
We want to show the following:
\begin{lemma}
\begin{enumerate}
\item 
$\phi^*(X_i)_{0\leq i\leq m}$ is an extension of  $(X'_j)_{0\leq j\leq m'}$.
\item
For $x'\in\supp(\phi^*(\cJ_i))$ we have the equalities of invariants \quad $\overline{\inv}(x')=\overline{\inv}(\phi_i(x'))$,\quad   $\nu(x')=\nu(\phi_i(x'))$ and $\rho(x')=\rho(\phi_i(x'))$.
\end{enumerate}
\end{lemma}
\noindent{\bf Proof.}
Denote by $s(\phi^*(X_i))$ the maximum number of $\phi^*(E)$ through a point in $\supp(\phi^*(\overline{\cJ_i}))$. 
In particular $s(\phi^*(X_i))\leq s_i$ for any index $0\leq i\leq m$.

Assume that for the index $l$ we can find index $l'$ such that $$\phi^*(X_{i_l},\overline{\cJ}_{i_l},E_{i_l},\overline{\mu})\simeq (X'_{j_{l'}},\overline{\cJ}'_{j_{l'}},E'_{j_{l'}},\overline{\mu}). \hfill{(\star)}$$
(This assumption is satisfied for  $l=0$.)
\begin{enumerate}

\item If $s(\phi^*(X_{i_l}))<s_{i_l}$ 
then the centers $C_i$ of blow-ups
 in the sequence $(X_i)_{i_l\leq i\leq i_{l+1}}$ are contained in the intersections of $s_{i_l}$ divisors in $E$ and do not hit the images  $\phi_i(\phi^*(X)_i)$. Thus $\phi^*(X_i)_{i_l\leq i\leq i_{l+1}}$ consists of isomorphisms. The property $(\star)$ will be satisfied for $l+1$ (and for the same $l'$).

\item If $s(\phi^*(X_{i_l}))=s_{i_l}>0$  then the intersections $(H')^s_{\alpha i}$ of $s=s(\phi^*(X_{i_l}))=s_{i_l}$ divisors 
are inverse images of $H^s_\alpha$ and the resolution process  
$\phi^*(X_i)_{i_l\leq i\leq i_{l+1}}$ is reduced to  resolving $\phi^*(\overline{\cJ}_{i_l|H^s_{\alpha i_l}})$.  
Moreover by the property $(\star)$,
 \begin{equation}\phi^*(\overline{\cJ}_{i_l|H^s_{\alpha i_l}})=\overline{\cJ}'_{j_{l'}|(H')^s_{\alpha j_{l'}}}\nonumber \end{equation}\noindent  for some $l'$ such that $s_{j_{l'}}=s_{i_l}$.

By commutativity of \'etale morphisms with the canonical resolution in lower 
 dimensions we know that all resolutions $\{(H')^s_{\alpha i}\}_{{i_l}\leq i \leq i_{l+1}}$ induced by  $\phi^*(X_i)_{i_l\leq i\leq i_{l+1}}$  are extensions of the canonical resolutions of
$\overline{\cJ}'_{i_{l'}|(H')^s_\alpha}$.
Moreover the restrictions $\phi_{i|(H')^s_{\alpha i}}: (H')^s_{\alpha i}\to H^s_{\alpha i}$ preserve the invariants $\inv$, $\mu$ and $\rho$.
 Thus $\phi^*(X_i)_{{i_l}\leq i \leq i_{l+1}}$ is an extension  of  $(X'_i)_{{j_{l'}}\leq i \leq j_{l'+1}}$. Moreover for ${i_l}\leq i <i_{l+1}$ and $x'\in \supp(\phi_i^*(\overline{\cJ}_i, \overline{\mu}))\cap (H')^s_{\alpha i}$, and $x=\phi_i(x')$ we have 
 \begin{equation}
  \overline{\inv}(x')\quad =\quad (s(x'),\, \inv_{({\phi_i^*(\overline{\cJ})}_{i|(H')^s_\alpha})}(x))\quad=\quad (s(x), \, \inv_{({\overline{\cJ}}_{i|H^s_\alpha})}(x))\quad=\quad \overline{\inv}(x).\nonumber\end{equation}
Analogously $\nu(x')=\nu(x)$ and $\rho(x')=\rho(x)$. The property $(\star)$ is satisfied for $l+1$ (and $l'+1$).

\item If $s(\phi^*(X_{i_k}))=s_{i_k}=0$ then 
the resolution process for $(X_i)_{i_k\leq i\leq m}$ is reduced 
to the canonical resolution of $\overline{\cJ}_{i_k|V(u)_{i_k}}$ on  a hypersurface of maximal contact $V(u)_{i_k}$.  Also the resolution process of $\overline{\cJ}'_{j_{k'}}\simeq \phi^*(\overline{\cJ}_{i_k})$
 is reduced to the
to canonical resolution of $\overline{\cJ}'_{i_k|V(u)_{i_k}}=\phi^*(\overline{\cJ}_{i_k|V(u)_{i_k}})$ on the hypersurface of maximal contact $V(u)_{i_k}$.
Since the inverse image of hypersurface of maximal contact is  a hypersurface of maximal contact by the same reasoning as before (replacing $H^s_\alpha$ with $V(u)$) we deduce that $\phi^*(X_i)_{i_k\leq i\leq m}$ is an extension of $(X'_i)_{i_{k'}\leq i\leq m}$. Moreover  for ${i_k}\leq i \leq m$ and $x'\in \supp(\overline{\cJ}'_i,\overline{\mu})$, 
$$\overline{\inv}(x')=\overline{\inv}(\phi_i(x')),\quad \nu(x')=\nu(\phi_i(x')),\quad\phi_i(\rho(x'))=\rho(\phi_i(x')).$$ 

\end{enumerate}
The lemma is proven.
\bigskip

{\bf Step 2}.   {\bf Resolving marked ideals $(X,{\cI},E,{\mu})$.}

For any  marked ideal $(X,{\cI},E,\mu)$ write 
 $$I=\cM({\cI}){\cN}({\cI}),$$  \noindent where $\cM({\cI})$ is the {\it monomial
part} of ${\cI}$, that is, the product of the principal ideals
defining the irreducible components of the divisors in $E$, and ${\cN}({\cI})$ is a {\it nonmonomial
part} which is not divisible by any ideal of a divisor in $E$.
Let $$\ord_{{\cN}({\cI})}:=\max\{\ord_x({\cN}({\cI}))\mid x\in\supp(\cI,\mu)\}.$$

\begin{definition} (Hironaka, Bierstone-Milman,Villamayor, Encinas-Hauser) 
By the {\it companion ideal} of $({\cI},\mu)$  where $I={\cN}({\cI})\cM({\cI})$ 
we mean the marked ideal of maximal order
\begin{displaymath}
O({\cI},\mu)= \left\{ \begin{array}{ll}
({\cN}({\cI}),\ord_{{\cN}({\cI})})\quad + \quad(\cM({\cI}),\mu-\ord_{{\cN}({\cI})})& \textrm{if   $\ord_{{\cN}({\cI})}<\mu$},\\ ({\cN}({\cI}),\ord_{{\cN}({\cI})})&\textrm{if   $\ord_{{\cN}({\cI})}\geq \mu$}. 
\end{array} \right.
\end{displaymath}
\end{definition} 

In particular $O({\cI},\mu)=({\cI},\mu)$ for ideals $({\cI},\mu)$ of maximal order.

{\bf Step 2a}. {\bf Reduction to the monomial case by using companion ideals.}

 By Step 1 we can resolve the marked ideal of maximal order $(\cJ,\mu_\cJ):=O({\cI},\mu)$ using the invariant $\overline{\inv}_{O({\cI},\mu)}$. By Lemma \ref{le: operations}, for any multiple test  blow-up of  $O({\cI},\mu)$, 
 
\centerline{$
\supp(O({\cI},\mu))_i= \supp[{\cN}({\cI}),\ord_{{\cN}({\cI})}]_i\,\, \cap \,\, \supp[M({\cI}),\mu-\ord_{{\cN}(H{\cI})}]_i =$}
 \centerline{$
 \supp[{\cN}({\cI}),\ord_{{\cN}({\cI})} ]_i\,\,\cap \,\, \supp({\cI}_i,\mu).$} 
 Consequently, such a resolution leads to the ideal $({\cI}_{r_1},\mu)$ such that $\ord_{{\cN}({\cI}_{r_1})}<\ord_{{\cN}({\cI})}$. This resolution is controlled by the invariants $\inv$, $\nu$ and $\rho$  defined for all $x\in  \supp({\cN}({\cI}),\ord_{{\cN}({\cI})})\cap \supp({\cI}_i,\mu)_i$, $$\inv(x)=\bigg(\frac{\ord_{{\cN}({\cI})}}{\mu},\overline{\inv}_{O({\cI},\mu)}(x)\bigg),\quad \nu(x)={\nu}_{O({\cI},\mu)}(x),\quad \rho(x)={\rho}_{O({\cI},\mu)}(x).$$ 
Then we repeat the procedure for $({\cI}_{r_1},\mu)$.  
 We find marked ideals $({\cI}_{r_0},\mu)=(\cI, \mu), ({\cI}_{r_1},\mu), \ldots, ({\cI}_{r_m},\mu)$ such that
$\ord_{{\cN}({\cI}_0)}>\ord_{{\cN}({\cI}_{r_1})}>\ldots>\ord_{{\cN}({\cI}_{r_m})}$.
  The procedure terminates after a finite number of steps when we arrive at the ideal $({\cI}_{r_m},\mu)$ with $\ord_{{\cN}({\cI}_{r_m})}=0$ or with $\supp({\cI}_{r_m},\mu)=\emptyset$. In the second case we get the resolution. In the first case $\cI_{r_m}=\cM({\cI}_{r_m})$ is monomial.

In Step 2a all points $x\in\supp(\cI,\mu)$ for which $\ord_c(\cI)\neq 0$ were assigned the invariants $\inv,\mu,\rho$. They are upper semicontinuous by the semicontinuouity of  $\ord_x$ and of the invariants $\overline{\inv},\mu,\rho$ for the marked ideals of maximal order.

{\bf Step 2b}. {\bf Monomial case $\cI=\cM({\cI})$}.

Define the invariants $$\inv(x)=(0,\ldots,0,\ldots),\quad \nu(x)=\frac{\ord_x({\cI})}{\mu}.$$

Let $x_1,\ldots,x_k$ define equations of the components ${D}^x_1,\ldots,{D}^x_k\in E$ through $x\in \supp(X,{\cI},E,\mu)$ and
 ${\cI}$  be generated by the monomial $x^{a_1,\ldots,a_k}$ at $x$. In particular $\nu(x)=\frac{a_1+\ldots+a_k}{\mu}$. 
 
 Let 
 $\rho(x)=\{D_{i_1},\ldots,D_{i_l}\}\in \Sub(E)$  be the maximal subset satisfying the properties

\begin{enumerate} 
\item  $a_{i_1}+\ldots+a_{i_l}\geq \mu.$
\item For any $j=1,\ldots,l$, $a_{i_1}+\ldots +\check{a}_{i_j}+\ldots +a_{i_l} < \mu. $ 
\end{enumerate}

Let $R(x)$ denote the subsets in $\Sub(E)$ satisfying the properties (1) and (2).
The maximal components of the $\supp({\cI},\mu)$ through $x$ are described by the intersections 
$\bigcap_{D\in A} D$ where $A\in R(x)$. The maximal locus of $\rho$ determines at most o one maximal component of $\supp({\cI},\mu)$ through each $x$. 
 
  After the blow-up at the maximal locus $C=\{x_{i_1}=\ldots=x_{i_l}=0\}$ of $\rho$, the ideal $\cI=(x^{a_1,\ldots,a_k})$ is equal to $\cI'=({x'}^{a_1,\ldots,a_{i_j-1},a,a_{i_j+1},\ldots,a_k})$ in the neighborhood  corresponding to $x_{i_j}$, where $a=a_{i_1}+\ldots+a_{i_l}-\mu<a_{i_j}$. In particular the invariant $\nu$ drops for all  points of some maximal components of $\supp({\cI},\mu)$. Thus the maximal value of $\nu$ on the maximal components of  $\supp({\cI},\mu)$ which were blown up is bigger than the maximal value of $\nu$ on the new maximal components of $\supp({\cI},\mu)$.
Since the set $\frac{1}{\mu}{\bf Z}_{\geq 0}$ of values of $\nu$ is discrete the algorithm 
terminates after a finite number of steps.\qed

\bigskip 
{\bf Commutativity of the resolution procedure in Step 2 with \'etale morphisms.}
The reasoning is the same as in Step 1.
Let $\phi: X'\to X$ be an \'etale morphism.  In Step 2a we find a sequence $r_0:=0<r_1<\ldots<r_k=r$ such that 
$\ord_{{\cN}({\cI}_{r_0})}>\ord_{{\cN}({\cI}_{r_1})}>\ldots>\ord_{\cN({\cI}_{r_m})}$
and for
 ${r_j}\leq i <r_{j+1}$,  $\ord_{{\cN}({\cI}_{r_j})}=\ord_{{\cN}({\cI}_{i})}$. Moreover the resolution process for $(\cI_i)_{r_j\leq i\leq r_{j+1}}$ is reduced to resolving the marked ideal of maximal order $O(\cI_{r_j})$.
Let $\ord_{{\cN}({\cI_{p_0}}')}>\ord_{{\cN}({\cI_{p_1}}')}>\ldots>\ord_{\cN({\cI_{p_{k'}}'})}$  be the  corresponding sequence defined for the canonical resolution of $({\cI'},{\mu})=\phi^*(\cI,\mu)$. 
\begin{lemma}
\begin{enumerate}
\item $(\phi^*(X_i))_{0\leq i\leq m}$ is an extension of $(X'_i)_{0\leq j\leq m'}$.
\item For $x'\in\supp(\phi^*(\cI_i))$ we have the equalities of invariants \quad ${\inv}(x')={\inv}(\phi_i(x'))$,\quad   $\nu(x')=\nu(\phi_i(x'))$ and $\rho(x')=\rho(\phi_i(x'))$.

\end{enumerate}
\end{lemma}
\noindent{\bf Proof.}
Note that all morphisms $\phi_i:\phi^*(X_i)\to X_i$ preserve the order of the nonmonomial part at a point $x\in \supp(\phi^*_i(\cI_i))$.
 Assume that for the index $l$ we can find an index $l'$ such that \\ \centerline{$\phi^*(X_{r_l},{\cI}_{r_l},E_{r_l},\overline{\mu})\simeq (X'_{p_{l'}},{\cI}'_{p_{l'}},E'_{p_{l'}},\overline{\mu})$.}

\begin{enumerate}
\item If $\ord_{\cN({\cI}_{r_{l}})}>\ord_{\cN(\phi_{r_l}^*{\cI}_{r_{l})}}=\ord_{\cN({\cI}_{p_{l'}})}$ then the centers of blow-ups of $(X_i)_{{r_l}\leq i <r_{l+1}}$ are contained in the loci of the points $x$ for which $\ord_x(\cN(\cI_i))=\ord_{\cN({\cI}_{r_{l})}}$. Therefore they are disjoint from images of $\phi^*(X_i)$. Consequently, $\phi^*(X_i)_{{r_l}\leq i <r_{l+1}}$ consists of isomorphisms.
\item If $\ord_{\cN({\cI}_{r_{l}})}=\ord_{\cN(\phi_{r_l}^*{\cI}_{r_{l}})}=\ord_{\cN({\cI}_{p_{l'}})}$ then $\phi_{r_l}^*(O({\cI}_{r_l}))=O(\phi_{r_l}^*({\cI}_{r_l})$.  
 By commutativity of the canonical resolution in Step 1 we get for any $x'\in\supp(O(\phi_i^*({\cI}_{i}))$,  and any ${r_l}\leq i <r_{l+1}$, $$\inv(x')\quad=\quad (\ord_{{\cN}(\phi^*{\cI}_{i}}), \, \overline{\inv}_{O(\phi^*({\cI}_{i}))}(x'))\quad=\quad
 (\ord_{\cN({\cI}_{i})}, \,  \overline{\inv}_{O({\cI}_{i})}(\phi_{i}(x')))\quad=\quad \inv(\phi_{i}(x')).$$
Analogously  $\nu(x')=\nu(\phi_i(x'))$, $\phi_i(\rho(x'))=\rho(\phi_i(x'))$ and $\phi^*(X_i)_{{r_l}\leq i <r_{l+1}}$ is an extension of the part of the  resolution of $(X'_i)_{{p_l}\leq i <p_{l+1}}$.
\item If $\cI_{r_k}=\cM(\cI_{r_k})$ and $\cI'_{p_{k'}}=\phi^*({\cI}_{r_k})=\cM(\phi^*({\cI}_{r_k}))$ are monomial the resolution process is controlled by the invariant $\rho$.  The set of values of $\rho$ on $X'$ can be identified via $\phi^*$
with subset of  a set of values of $\rho$ on $X$: $\Sub(E')\subset \Sub(E)$.
By definition $\rho$ and $\nu$ commute with smooth morphisms: $\rho(\phi(x'))=\rho(x')$ and $\nu(\phi(x'))=\nu(x')$.
The blow-ups on $(X)_{r_{k}\leq i \leq m}$ are performed at the centers where $\rho$ attains its maximum. Thus the induced morphisms on  $\phi^*(X)_{r_{k}\leq i \leq m}$ either are blow-ups performed at the centers where $\rho$ attains a maximum or are isomorphisms.
Consequently, $\phi^*(X_i)_{r_{k}\leq i \leq m}$  is an extension of $(X'_i)_{p_{k'}\leq i \leq m'}$.  \qed
\end{enumerate}

\bigskip
\subsection{Summary of the resolution algorithm} 

The resolution algorithm can be represented by the following scheme.

\noindent {\bf Step 2}.  Resolve $(\cI,\mu)$.

{\bf Step 2a.} Reduce $(\cI,\mu)$ to the monomial marked 
 ideal $\cI=\cM(\cI)$.
$$\Downarrow$$
\quad\quad If $\cI\neq\cM(\cI)$, decrease the maximal order of the nonmonomial part $\cN(\cI)$ by resolving the companion 

\quad\quad ideal $O(\cI,\mu)$. For $x\in\supp(O(\cI,\mu))$, set
$$\inv(x)=(\ord_x(\cN(\cJ))/\mu,\overline{\inv}_{O(\cI,\mu)}).$$

\quad\quad{\bf Step 1}. Resolve the companion ideal $(\cJ,\mu_\cJ):=O(\cI,\mu)$ : 

\quad\quad\quad  Replace $\cJ$ with $\overline{\cJ}:=\cC(\cH(\cJ))\simeq \cJ$. (*)

\quad\quad\quad{\bf Step 1a}. Move apart all strict transforms of $E$ and $\supp(\overline{\cJ},\mu)$.
$$\Downarrow$$
\quad\quad\quad\quad\quad   Move apart all intersections $H^s_\alpha$
of $s$ divisors in $E$

\quad\quad\quad\quad\quad
(where $s$ is the maximal number of divisors in $E$ through 
 points in $\supp(\cI,\mu)$).
 $$\Updownarrow$$

\quad\quad\quad\quad\quad{\bf Step 1aa}. If $\overline{\cJ}_{|H^s_\alpha}=0$ for some $H^s_\alpha\subset \supp(\overline{\cJ})$, blow up
$H^s_\alpha$.  For  $x\in  H^s_\alpha$ set $$\inv(x)=(\ord_x(\cN(\cJ))/\mu,s,\infty,0,\ldots),\quad \nu(x)=0,\quad \rho(x)=\emptyset.$$

\quad\quad\quad\quad\quad{\bf Step 1ab.} If $\overline{\cJ}_{|H^s_\alpha}\neq 0$ for any $\alpha$,  resolve all $\overline{\cJ}_{|H^s_\alpha}$. For  $x\in \supp(\overline{\cJ},\mu)\cap H^s_\alpha$ set
$$\quad \inv(x)=(\ord_x(\cN(\cJ))/\mu,s,\inv_{\overline{\cJ}_{|H^s_\alpha}}(x)),\quad \nu(x)=\nu_{\overline{\cJ}_{|H^s_\alpha}}(x),\quad \rho(x)=\rho_{\overline{\cJ}_{|H^s_\alpha}}(x).$$

\quad\quad\quad\quad\quad Blow up the centers where $(\inv,\rho)$ is maximal. 

\quad\quad\quad{\bf Step 1b} If  the strict transforms of $E$ do not intersect  $\supp(\overline{\cJ},\mu)$, resolve $(\overline{\cJ},\mu)$.
$$\Updownarrow$$

\quad\quad\quad\quad\quad{\bf Step 1ba.} If the set $\codim(1)(\supp(\overline{\cJ})$ of   codimension one components  is nonempty, 

\quad\quad\quad\quad\quad\quad blow it up.
For $x\in \supp(\overline{\cJ},\mu)=\codim(1)(\supp(\overline{\cJ})$ set 
$$\quad \inv(x)=(\ord_x(\cN(\cJ))/\mu,0,\infty,0),\quad \nu(x)=0,\quad \rho(x)=\emptyset.$$

\quad\quad\quad\quad\quad{\bf Step 1bb.} Simultaneously resolve
 all $\overline{\cJ}_{|V(u)}$ (by induction), where $V(u)$ is a hypersurface of 

\quad\quad\quad\quad\quad\quad maximal contact. For $x\in \supp(\overline{\cJ},\mu)\setminus \codim(1)(\supp(\overline{\cJ}))$ 
set
$$\quad\quad\quad \inv(x)=(\ord_x(\cN(\cJ))/\mu,s(x),\inv_{\overline{\cJ}_{|V(u)}}(x)),\quad \nu(x)=\nu_{\overline{\cJ}_{|V(u)}}(x),\quad \rho(x)=\rho_{\overline{\cJ}_{|V(u)}}(x).$$

 \quad\quad\quad\quad\quad\quad  Blow up the centers where $(\inv,\rho)$ is maximal. 


{\bf Step 2b}. Resolve the monomial marked ideal $\cI=\cM(\cI)$.

(Construct the invariants $\inv$, $\rho$ and $\nu$ directly for $\cM(\cI)$.)

\begin{remarks}

\bigskip
\begin{enumerate} 

\item (*) The ideal $\cJ$ is replaced with $\cH(\cJ)$ to ensure that the invariant constructed in Step 1b is independent of the choice of the tangent direction $u$.

We replace $\cH(\cJ)$  with $\cC(\cH(\cJ))$ to ensure the equalities $\supp(\cJ_{|S})=\supp(\cJ)\cap S$, where $S=H^s_{\alpha}$ in Step 1a and  $S=V(u)$ in Step 1b. 
\item If $\mu=1$ the companion ideal is equal to $O(\cI,1)=(\cN(\cI),\mu_{\cN(\cI)})$ so the general strategy of the resolution of $\cI,\mu$ is to decrease the order of the nonmonomial part and then to resolve the monomial part.
\item In particular if we desingularize $Y$ we put  $\mu=1$ and $\cI=\cI_Y$ to be equal to the sheaf of the subvariety $Y$ and we resolve the marked
ideal $(X,\cI,\emptyset,\mu)$. The nonmonomial part $\cN(\cI_i)$ is nothing but the \underline{weak transform} $(\sigma^i)^{\rm w}(\cI)$ of $\cI$. 
\end{enumerate}
\end{remarks}

\bigskip 
\subsection{Desingularization of plane curves}
We  briefly illustrate the resolution procedure for plane curves.

Let $C\subset {\bf A}^2$ be a plane curve defined by $F(x,y)=0$ (for instance $x^2+y^5=0$).
We assign to the curve $C$ the marked ideal $(X,\cI_C,\emptyset,1)$.
The nonmonomial part of a controlled transform of the ideal $\cI_C$ is the ideal of the strict transform of the curve (In general it is the weak transform of the subvariety). In particular $\cI_C=\cN(\cI_C)$. 

In {\bf Step 2a} we  form the companion ideal which is equal to $\cJ:=O(\cI_C)=(\cI_C,\mu)$, where $\mu$ is the maximal multiplicity. Resolving $O(\cI_C)$ will eliminate the maximal multiplicity locus of $C$ and  
decrease the maximal multiplicity of the ideal of the strict transform of $C$.
 The maximal multiplicity locus of $C$ 
is defined by $\supp(\cI_{C},\mu)=V(\cD^{\mu-1}(\cI_C))$,  which is a finite set of points for a singular curve. 

In the example $\mu=2$ and  $\cJ=O(\cI_C)=(\cI_C,2)$, $T(\cJ)=(\cD(\cI_C),1)=((x,y^4),1)$, $\supp(\cI_{C},2)=V(x,y^4)=\{(0,0)\}$.

In {\bf Step 1} we resolve the companion ideal $\cJ=(\cI_C,\mu)$. 
First   replace $\cJ$ with $\overline{\cJ}:=\cC(\cH({\cJ}))$. In the example $$\overline{\cJ}:=\cC(\cH(\cJ))=\cH({\cJ})=(x^2,xy^4,y^5).$$
Since at the beginning there are no exceptional divisors we go directly to Step 1b.

{\bf Step 1b}. For any point $p$ with multiplicity $\mu$ we find a tangent direction $u\in T(\cI,\mu)$ at $p$. In particular $u=x$ for $p=(0,0)$. Then assign to $p$ the invariant
$$\inv(p)=(\mu,0,\ord_p(\cJ_{|V(u)})/\mu!,\infty,0,\ldots).$$
In general for local coordinates $u,v$ at $p$ we have $\cJ_{|V(u)}=(v^m,\mu!)$ and we can write the invariant as 
$$\inv(p)=(\mu,0,m/\mu!,\infty,0,\ldots),\quad \nu(p)=0,\quad \rho(p)=\emptyset.$$
In the example $\cJ_{|V(u)}=\cJ_{|V(x)}=(y^5,2)$ and
$$\inv(p)=(2,0,5/2,\infty,0,\ldots).$$
The resolution of $\cJ_{|V(u)}$ consists of two steps:
Reducing to the monomial case in {\bf Step 2a} and resolving the monomial case in {\bf Step 2b}.
 We blow up all points for which this invariant is maximal.
After the blow-ups $\cJ_{|V(u)}$ is transformed as follows:
$$(v^m,\mu!)\mapsto (y_{\rm exc}^{m-\mu!},\mu!)$$
If $\supp(y_{\rm exc}^{m-\mu!},\mu!)=\emptyset$ then $\cJ_{|V(u)}$ is resolved and the multiplicity of the corresponding points drops. Otherwise  $\sigma^{\rm c}(\cJ)_{|V(u)}=(y_{\rm exc}^{m-\mu!},\mu!)$ is monomial for all points with  the highest multiplicity. The assigned invariant is 
$$\inv(p')=(\mu,0,0,0,0,\ldots),\quad \nu(p')=(m-\mu!)/\mu!, \quad \rho(p')=D_{\rm exc}.$$
In the example $\sigma^{\rm c}(y^5,2)=(y_{\rm exc}^3,2)$ and
$\inv(p')=(2,0,0,\ldots)$ and $\nu(p')=3/2$.
The equation of the strict transform of $C$ at the point with the highest multiplicity changes as follows
\begin{equation} x^2+y^5=0 \mapsto x^2+y^3_{\rm exc}=0 \end{equation}
After the next blow-up the invariants for all points with the highest multiplicity are equal
$$\inv(p'')=(\mu,0,0,0,0,\ldots)\quad \nu(p')=(m-2\mu!)/\mu! \quad \rho(p'')=D'_{\rm exc}$$
We continue blow-ups till  $m-l\mu!\leq\mu!$. At this moment  $\supp(\sigma^{\rm c}(\cJ)_{|V(u)})=\emptyset$  and the marked ideal $\cJ_{|V(u)}$ is resolved (as in Step 1b). Resolving  $\cJ_{|V(u)}$ is equivalent to resolving $\cJ$ and results in dropping the maximal multiplicity. In the example after the second blow-up
$5-2\cdot 2\leq 2$ and the maximal multiplicity drops to $1$. 
\begin{equation} x^2+y^3_{\rm exc}=0 \mapsto x^2+y'_{\rm exc}=0 \end{equation}
After all points with the highest multiplicity are eliminated and the maximal multiplicity of points drops we reconstruct our companion ideals for the controlled transform of $\cI_C$.
The companion ideal of $\sigma^{\rm c}(\cI_C)$ is equal to $\cJ':=(\cI_{C'},\mu')$, where $\cI_{C'}$ is the ideal of the strict transform and $\mu'$ is the highest multiplicity. As before $\supp(\cJ')$ defines the set of  points with the highest multiplicity. In our example the curve $C'$ is already smooth and $\mu'=1$. However the process of the embedded desingularization is not finished at this stage. Some exceptional divisors may pass through  the points with the highest multiplicity.  In the course of resolution of $\cJ'$ we first move apart all strict transforms of the exceptional divisors and the set of points with multiplicity $\mu'$.
This is handled in {\bf Step 1a} by resolving $\overline{\cJ'}_{|H^s_\alpha}$. The maximum number of the exceptional divisors passing through points of $\supp(\cJ')$ can be $s=2$ or $s=1$.
If $s=2$ then the assigned invariants are
$$\inv(p')=(\mu',2,\infty,0,\ldots),\quad \nu(p)=0,\quad \rho(p)=\emptyset.$$
The blow-up of the point separates the divisors.
If $s=1$ then $H^s_\alpha=D_\alpha$ is a single divisor,
$$\inv(p')=(\mu',1,\ord_{p'}(\cN(\overline{\cJ}'_{|D_\alpha}))/\mu!,0,\ldots) ,\quad \nu(p')=0,\quad \rho(p')=\emptyset,$$
\noindent where $\cN(\overline{\cJ}'_{|D_\alpha})=((v^m),\mu'!)$. We resolve  this ideal as above: First we eliminate the nonmonomial part $\cN(\overline{\cJ}'_{|D_\alpha})$ and then resolve the resulting monomial ideal $\overline{\cJ}'_{|D_\alpha}$.

In our example the second exceptional divisor $y'_{\rm exc}=0$ passes through the point $p''$: $x=y'_{\rm exc}=0$.  
\begin{equation}\cN(\overline{\cJ}'_{|D_\alpha})=(x^2,1)\quad \mapsto \quad   (y''_{\rm exc},1)\quad \mapsto \quad (\cO_{D'},1)\nonumber\end{equation}
\begin{equation} x^2+y'_{\rm exc}=0\quad  \mapsto \quad y''_{\rm exc}+y'_{\rm exc}=0 \quad \mapsto \quad 1+y'_{\rm exc}=0 \nonumber\end{equation}
$$\inv(p'')=(1,1,2,0\ldots), \nu(p'')=0,\quad  \mapsto \quad \inv(p''')=(1,1,0,\ldots), \nu(p''')=1$$

After the ideals are resolved the strict transforms of all exceptional diviors  are moved away from the set of points with highest multiplicity and we arrive at Step 1b.
If $\mu'=1$ we stop the resolution procedure. At this moment the invariant for all points of the strict transform of $C$ is constant and equal to
$\inv(p)=(1,0,\infty,0,\ldots),\quad \mu(p)=0$.
The strict transform of $C$ is now smooth and has simple normal crossings with exceptional divisors. It defines a hypersurface of maximal contact.  

If $\mu'>1$ we repeat the procedure for Step 1b described above. After the ideals $\cJ_{|V(u')}$ are resolved the highest multiplicity drops. The procedure terminates when the invariant is constant along $C$ and equal  to
$$\inv(p)=(1,0,\infty,0,\ldots),\quad\quad \mu(p)=0.$$

\section{Conclusion of the resolution algorithm} 
\subsection{Commutativity of resolving  marked ideals with smooth morphisms} 
Let $(X,\cI,\emptyset,\mu)$ be a marked ideal and $\phi: X'\rightarrow X$ be a smooth morphism of relative dimension $n$. Since the canonical resolution  is defined by the invariant it suffices to show that
$\inv(\phi(x))=\inv(x)$. Let $U'\subset X'$ be a neighborhood of $x$ such that there is afactorization
$\phi: U'\buildrel \phi'\over
\to X\times {\bf A}^n \buildrel  \pi \over
\to X$, where $\phi'$ is \'etale and $\pi$ is the natural projection. 
The canonical resolution $(X_i\times {\bf A}^n$) of $p^*(X,\cI,E,\mu)$ is induced by the canonical resolution $(X_i)$ of $(X,\cI,E,\mu)$ and the invariants $\inv$, $\mu$ and $\rho$ are preserved by $\pi$.
Then for $x'\in \supp(X',\cI',E',\mu)\cap U'$ we have $\inv(\phi(x'))=\inv(\pi(\phi' (x')))=\inv(\phi' (x'))$. Since
$\phi'$ is \'etale,  the resolution $\phi_{|U'}^*(X_i)={\phi'}^*(X_i\times {\bf A}^n)$ is an extension of the canonical resolution of $\cI'_{|U'}$ and
$\inv(\phi' (x))=\inv(x)$. Finally $\inv(\phi(x))=\inv(x)$. Analogously $\mu(\phi(x))=\mu(x)$ and $\rho(\phi(x))=\rho(x)$.

\subsection{Commutativity of resolving marked ideals $(X,\cI,\emptyset,1)$ with embeddings of ambient varieties}

Let $(X,\cI,\emptyset,1)$ be a marked ideal and $\phi: X\hookrightarrow X'$ be a closed embedding of smooth varieties. Then $\phi$ defines the marked ideal $(X',\cI',\emptyset,1)$, where $\cI'=\phi_*(\cI)\cdot\cO_{X'}$. We may assume that $X$ is a subvariety of $X'$ locally generated by parameters $u_1,\ldots, u_k$. Then $u_1,\ldots, u_k\in\cI'(U')=T(U')$ define tangent directions on some open $U'\subset X'$.  We run steps 2a and 1bb of our algorithm
 through. In step 2a we  assign $\inv(x)=(1,\overline{\inv}_\cI(x))$ (since the maximal order of ${\cI}=\cN(\cI)$ is equal to $1$, and $\cI=O({\cI})$) and in Step 1bb (nonboundary case $s(x)=0$) we assign
\begin{equation}
\inv_{\cI'}(x)=(1,0,{\inv}_{\cI'_{|V(u_1)}}(x)),\quad \nu_{\cI'}(x)=\nu_{\cI'_{|V(u_1)}}(x) ,\quad \rho_{\cI'}(x)=\rho_{\cI'_{|V(u_1)}}(x)
\end{equation}  
passing to the hypersurface $V(u_1)$.
By Step 1bb resolving $(X',\cI',\emptyset,\mu)$ is locally equivalent to resolving $(V(u),\cI'_{V(u)},\emptyset,\mu)$ with relation between invariants defined by (1).
By repeating the procedure $n$ times and restricting to the tangent directions $u_1,\ldots, u_k$ of the marked ideal $\cI$ on $X$ we obtain: 
\begin{equation}\inv_{\cI'}(x)=(1,0,1,0,\ldots,1,0,\inv_{\cI}(x)),\quad \rho_{\cI'}(x)=\rho_{\cI}(x),\quad \nu_{\cI'}(x)=\nu_{\cI}(x).\end{equation}
Resolving  $(X',\cI',\emptyset,\mu)$ is  equivalent to resolving $(X,\cI,\emptyset,\mu)$ with relation between invariants defined by (10).

\bigskip

\subsection{Commutativity of resolving   marked ideals  with isomorphisms not preserving  the ground field}
\begin{lemma} Let $X,X'$ be varieties over $K$ and $\cI$ be a sheaf of ideals on $X$. Let  $\phi: X'\to X$ be an isomorphism  over ${\bf Q}$. Then $\phi^*(\cD^i(\cI))=\cD^i(\phi^*(\cI))$ for any $i$. 
\end{lemma} 
\noindent{\bf Proof.} It suffices to consider the case $i=1$. The sheaf $\cD^i(\cI)$ is locally
generated by functions $f\in \cI$ regular on some open subsets $U$ and their derivatives $D(f)$ . Then $\phi^*(\cD^i(\cI))$ is locally 
generated by $\phi^*(f)$ and $\phi^*Df=\phi^*D(\phi^*)^{-1}\phi^*f$. But for any derivarition $D\in \Der_k(\cO(U))$, $D':=\phi^*D(\phi^*)^{-1}\in \Der_k(\cO(\phi^{-1}(U))$ defines a $K$-derivation of
$\cO(\phi^{-1}(U))$.\qed

\begin{proposition} Let $(X,\cI,E,\mu)$ be a marked ideal. Let $\phi: X'\to X$ be an isomorphism over ${\bf Q}$. For any canonical resolution $(X_i)$ of $X$ the induced resolution $(X_i'):=(X_i\times_X X')$ is canonical. 
Moreover the isomorphism $\phi$ lifts to isomorphisms $\phi_i:X'_i\to X_i$ such that
$$\inv(\phi(x))=\inv(x),\quad \nu(\phi(x))=\nu(x), \quad \rho(\phi(x))=\phi(\rho(x)).$$
\end{proposition}

\noindent{\bf Proof.} Induction on dimension of $X$. First assume that $(X,\cI,E,\mu)=(X,\cJ,E,\mu)$ is of maximal order as in Step 1. Then, by the lemma, $\phi^(\cC(\cH(\cJ)))=\cC(\cH(\cJ))$.
 Resolution algorithm in Step 1a is reduced to resolution of the restrictions of marked ideals $\cI$ to intersections of the exceptional divisors $H^s_\alpha$. This procedure commutes with the isomorphism $\phi$. Moreover $$\overline{\inv}(\phi(x)\quad =\quad(s(\phi(x)), {\inv}_{|\phi(H^s)}(\phi(x)))\quad =\quad (s(x),\inv_{|H^s}(x))\quad =\quad\overline\inv(x),$$ $$\ \nu(\phi(x))\quad =\quad\nu_{|\phi(H^s)}(\phi(x))\quad =\quad\nu_{|H^s}(x)\quad =\quad\nu(x), \quad \rho(\phi(x))=\phi(\rho(x)),$$ by the inductive assumption.
In Step 1b we reduce resolution of the marked ideal to its restriction to a hypersurface of maximal contact defined by $u\in \cD^{\mu-1}(\cI)$ on an open subset $U$. This procedure commutes with $\phi$. The corresponding marked ideal $(\phi^*(\cJ),\mu)$ is restricted to the  hypersurface of maximal contact on $\phi^{-1}(U)$ defined by $\phi^*(u)\in \phi^*(\cD^{\mu-1}(\cJ))$.
The invariants defined in this step commute with $\phi$ by the inductive assumption.
 In Step 2a we decompose arbitrary marked ideal into the monomial and nonmonomial part. Since an isomorphism $\phi: X'\to X$ maps  divisors in $E'$ to  divisors in  $E$ it preserves this decomposition.
Consequently, it preserves companion ideals and the invariants $\inv,\rho,\mu$ defined in  Step 2a. Also the invariants defined in  Step 2b are preserved by $\phi$. Therefore $\phi$ commutes with canonical resolutions.

\subsection{Resolving marked ideals over  a non-algebraically closed field}
Let $(X,\cI,E,\mu)$ be a marked ideal defined over a field $K$. Let $\overline{K}$ be the algebraic closure of $K$.  Then the base change $K\mapsto\overline{K}$ defines the $G=\Ga(\overline{K}/K)$-invariant marked ideal $(\overline{X},\overline{\cI},\overline{E},\mu)$ (over $\overline{K}$). The canonical resolution $(\overline{X},\overline{\cI},\overline{E},\mu)$ is $G$-equivariant and defines canonical resolution of $(X,\cI,E,\mu)$ over $K$. This resolution commutes with smooth morphisms and embeddings of the ambient varieties over $K$ and with isomorphisms over $\bf Q$. 

\subsection{ Principalization} 

Resolving  the marked ideal $(X, {\cI},\emptyset,1)$ determines a principalization commuting with smooth morphisms, group actions and embeddings of the ambient varieties. 
 
 The principalization is often reached at an earlier stage upon transformation to the monomial case (Step 2b). This moment is detected by the invariant $\inv$, which becomes equal to $\inv(x)=(0,\ldots,0,\ldots).$ (However the latter procedure does not commute with embeddings of ambient varieties)

\subsection{Weak embedded desingularization}

Let $Y$ be a closed subvariety of the variety $X$. Consider the marked ideal $(X, {\cI}_Y,\emptyset,1)$. Its support $\supp({\cI}_Y,1)$ is equal to $Y$. In the resolution process of $(X, {\cI}_Y,\emptyset,1)$, the strict transform of $Y$ is blown up. Otherwise  the generic points would be transformed isomorphically, which contradicts  the resolution of $(X, {\cI}_Y,\emptyset,1)$. At the moment where the strict transform is blown up the invariant along it is the same for all its points and 
equal to $$\inv(x)=(1,0;1,0;\ldots;1,0;\infty;0,0,0,\ldots),$$
\noindent where $(1,0)$ is repeated $n$ times. This value of invariant can be computed for the generic smooth point of $Y$. We apply  Step 2a ($ \ord_x({\cI})=1, \cI=O({\cI})$) and Step 1b (nonboundary case $s(x)=0$) passing to a hypersurface $n$ times. 
Each time after running through 2a and 1b we adjoin a couple (1,0) to the constructed invariant. 
After the n-th time we arrive at Step 1ba where the algorithm terminates and $\infty$ followed by zeros is added at the end of the invariant.

\bigskip
\subsection{Bravo-Villamayor strengthening of the Weak Embedded Desingularization}

\begin{theorem}(Bravo-Villamayor (see \cite{BV}, \cite{BM3})) Let  $Y$  be a reduced closed subscheme of a smooth variety $X$ and  $Y=\bigcup Y_i$ be its decomposition into the union of irreducible components.
There is a canonical resolution of a subscheme $Y\subset X$, that is, a sequence of blow-ups $(X_i)_{0\leq i\leq r}$ subject to conditions (a)-(d) from Theorem \ref{th: 2} such that the strict transforms $\widetilde{Y}_i$ of $Y_i$ are smooth and disjoint. Moreover the full transform of $Y$ is of the form $$(\widetilde{\sigma})^*(\cI_Y)= \cM((\widetilde{\sigma})^*(\cI_Y))\cdot\cI_{\widetilde{Y}},$$ where $\widetilde{Y}:=\bigcup\widetilde{Y}_i\subset\widetilde{X}$ is a disjoint union of the strict transforms $\widetilde{Y}_i$ of $Y_i$, $\cI_{\widetilde{Y}}$ is the sheaf of ideals of $\widetilde{Y}$ and $\cM((\widetilde{\sigma})^*(\cI_Y))$ is the monomial part of $(\widetilde{\sigma})^*(\cI_Y)$.
\end{theorem}

\noindent{\bf Proof.} 
Consider the canonical resolution procedure for the marked ideal $(X,\cI_Y,\emptyset,1)$ (and in general for $(X,\cI,E,\mu)$) described in the proof of Proposition \ref{pr: resol}.
We shall modify the construction of the invariants in the canonical resolution.
In Step 1 we define $\overline{\inv'}$, $\rho'$, $\nu'$ in the same way as before.
In Step 2 we modify the definition of the companion ideal to be

\begin{displaymath}
O'({\cI},\mu)= \left\{ \begin{array}{ll}({\cM}({\cI}),1) &\textrm{if   $\ord_{{\cN}({\cI})}\leq 1$ and $\mu=1$ and ${\cM}({\cI})\neq \cO_X$},\\ 
O({\cI},\mu)
& \textrm{otherwise}.\\ 

\end{array} \right.
\end{displaymath}

We define invariants as follows. 
If   $\ord_{{\cN}({\cI})}\leq 1$ and $\mu=1$ and $({\cM}({\cI}))\neq \cO_X$ we set 
\begin{equation} \inv'(x)=(3/2,0,0,\ldots),\quad \nu'(x)=\nu_{(\cM(\cI),1)}(x),\quad \rho'(x)=\rho_{(\cM(\cI),1)}(x)\nonumber
\end{equation}
\noindent defined as in Step 2b. Otherwise we put as before
\begin{equation} \inv(x)=\bigg(\frac{\ord_{{\cN}({\cI})}}{\mu},\overline{\inv}_{O({\cI},\mu)}(x)\bigg),\quad \nu(x)={\nu}_{O({\cI},\mu)}(x),\quad \rho(x)={\rho}_{O({\cI},\mu)}(x).\nonumber
\end{equation}

Note that the reasoning is almost the same as before. 
The difference occurs for resolving marked ideals $(\cI,1)$  in Step 2 when we arrive at the moment when $\max\{\ord_x(\cN(\cI))\}=1$. Let $\nu^1(x)$ denote the first coordinate of the invariant $\inv$. 
Note that $\nu^1(x)=3/2$ for all points of $\supp(\cI,1)$ for which $\cI$ is not purely nonmonomial.  
First resolve  its monomial part as in Step 2b (for all points with $\nu^1(x)=3/2$). The blow-ups are performed at exceptional divisors for which $\rho(x)$ is maximal. We arrive at the purely nonmonomial  case  ($\nu^1(x)=1$) and  continue resolution as before.

 Let us order the codimensions of the components $Y_i$ in a increasing sequence $r_1:=\codim{Y_1}\leq\ldots\leq r_k:=\codim{Y_k}$.

 We shall run  Steps 1-2 of this procedure with the above modifications 
 till the strict transform of one of the components $Y_i$ is the center of the next blow-up. At this point the invariants are constant along this strict transform and are equal to $$\inv(x)=(1,0;1,0;\ldots;1,0;\infty;0,0,0,\ldots),\quad\nu'(x)=0\quad\rho'(x)=\emptyset$$
\noindent where $(1,0)$ is repeated $r_1$ times. (These are the values of the invariants for a generic smooth point of $\widetilde{Y}_1$.) 

{\bf Claim}: {\it Let $(\cI,1)$ be a marked ideal on $X$ such that $Y_1$ is an irreducible component of $\supp(\cI)$. Moreover assume that $\cI=\cI_{Y_1}$ in a neighborhood of a generic point of $Y_1$. At the moment of the (modified resolution process) for which  $$\max(\inv(x))=(1,0;1,0;\ldots;1,0;\infty;0,0,0,\ldots)$$ \noindent  the controlled transform of $(\cI,1)$ is equal to $\cI_{\widetilde{Y}_1}$ in the neighborhood of the strict transform  ${\widetilde{Y}_1}$ of $Y_1$.} 

We prove this claim by induction on codimension.
Note that when we run Step 2a of the algorithm  at some point we arrive at a marked ideal for which $\max(\nu^1(x))=1$. At this stage
$\cI=\cN(\cI)$ is purely nonmonomial and $O'(\cI)=(\cI,1)$. Note that starting from this point the controlled transform of $\cI$ remains nonmonomial for all points with ${\nu^1(x)}=1$. 
Then we go to Step 1 and construct $\cC(\cH(\cI,1))=(\cI,1)$. In Step 1a we run the algorithm arriving at the nonboundary case in Step 1b. 
At this point we restrict $\cI$ to a smooth hypersurface of maximal contact $V(u)$. If we are
in Step 1ba this hypersurface is the strict transform of $Y_1$. Moreover the order of the controlled transform $\sigma^{\rm c}(\cI)$ of $\cI$ is  1 along the strict transform of $Y_1$ and thus $\sigma^{\rm c}(\cI)$ is the ideal of this strict transform (in the neighborhood of the strict transform). 

In case 1bb we apply the (modified) canonical resolution to the restriction $\cI_{|V(u)}$.
This restriction $\cI_{|V(u)}$ satisfies the assumption of the claim for $Y_1\subset V(u)$ (we
 skip indices here). 
 By the inductive assumption the controlled transform $\sigma^{\rm c}(\cI)_{|V(u)}$ of $(\cI_{|V(u)},1)$ is locally equal to the ideal  of the strict transform $\widetilde{Y}_1\subset \widetilde{V(u)}$ of $Y_1$. Since  $u\in \sigma^{\rm c}(\cI)$ it follows that $\sigma^{\rm c}(\cI)=\cI_{\widetilde{Y}_1}$
(in  the neighborhood of $\widetilde{Y}_1$). The claim is proven.

All the strict transforms of codimension $r_1$ are isolated. We  continue the (modified) canonical resolution procedure ignoring these isolated components. We arrive at the moment where some codimension $r_2>r_1$ component is the center of the blow-up and the invariant $\inv$ is equal to $$\inv(x)=(1,0;1,0;\ldots;1,0;\infty;0,0,0,\ldots),\quad\nu'(x)=0,\quad\rho'(x)=\emptyset,$$
\noindent where $(1,0)$ is repeated $r_2$ times.  Again by the claim the controlled transform of all codimension $r_2>r_1$ components coincide with the strict transform and are isolated. Starting from this moment those components are ignored in the resolution process. We continue for all $r_i$.
At the end we principalize all components if there are any which do not intersect the strict transforms of $Y_i$.
\bigskip

\subsection{Desingularization} 

Let $Y$ be an algebraic variety over $K$. 
By the compactness of $Y$ we find a cover of affine subsets $U_i$ of $Y$ such that each $U_i$ is embedded in an affine space ${\bf A}^{n}$ for $n>>0$. We can assume that the dimension
$n$  is the same for all $U_i$  by taking if necessary embeddings of affine spaces ${\bf A}^{k_i}\subset {\bf A}^n$. 
\begin{lemma} \label{le: l} Let $\phi_1,\phi_2: Y\subset {\bf A}^{n}$ be two embeddings defined by two sets of generators $g_1,\ldots,g_n$ and $h_1,\ldots,h_n$ respectively. Define three embeddings $\Psi_i:Y\to {\bf A}^{2n}$ for $i=0,1,2$ such that
\begin{eqnarray}
\nonumber \Psi_0(x)=(g_1(x),\ldots,g_n(x),h_1(x),\ldots,h_n(x)),\\
\nonumber \Psi_1(x)=(g_1(x),\ldots,g_n(x),0,\ldots,0),\\
\nonumber \Psi_2(x)=(0,\ldots,0,h_1(x),\ldots,h_n(x)).
\end{eqnarray}
Then there are automorphisms $\Phi_1,\Phi_2$ of ${\bf A}^{2n}$ such that for $i=1,2$,
$$\Phi_i\Psi_0=\Psi_i.$$
\end{lemma}
\noindent{\bf Proof.} Fix coordinates $x_1,\ldots,x_n,y_1,\ldots,y_n$ on ${\bf A}^{2n}$.
Find polynomials $w_i(x_1,\ldots,x_n)$ and $v_i(x_1,\ldots,x_n)$ such that for $i=1,\ldots, n$,
 \begin{equation} w_i(h_1,\ldots,h_n)=g_i, \quad  v_i(g_1,\ldots,g_n)=h_i. \nonumber\end{equation}
Set \begin{eqnarray} \Phi_1(x_1,\ldots,x_n,y_1,\ldots,y_n)=(x_1,\ldots,x_n,y_1-v_i(x_1,\ldots,x_n),\ldots,y_n-v_n(x_1,\ldots,x_n)), \nonumber \\ \Phi_2(x_1,\ldots,x_n,y_1,\ldots,y_n)=(x_1-w_1(y_1,\ldots,y_n),\ldots,x_n-w_n(y_1,\ldots,y_n),y_1,\ldots,y_n). \nonumber
\end{eqnarray}\qed
\begin{proposition} For any affine variety $U$ there is a smooth variety $\widetilde{U}$ along with a birational morphism $\res: \widetilde{U}\to U$ subject to the conditions:
\begin{enumerate}
\item For any closed embedding $U\subset X$ into a smooth affine variety $X$, there is an embedding $\widetilde{U}\subset \widetilde{X}$ into a smooth variety $\widetilde{X}$ which is a canonical embedded desingularization of $U\subset X$.
\item For any open embedding $V\hookrightarrow U$ there is an open embedding of resolutions $\widetilde{V}\hookrightarrow \widetilde{U}$ which is a lifting of $V\to U$ such that $\widetilde{V} \to \res_U^{-1}(V)$ is an isomorphism over $V$.

\end{enumerate}
\end{proposition}

\noindent{\bf Proof} 
(1) Consider a closed embedding of $U$ into a smooth affine variety $X$ (for example $X={\bf A}^n$). The canonical embedded desingularization $\widetilde{U}\subset \widetilde{X}$ of $U\subset X$ defines 
the desingularization $\widetilde{U}\to U$. This desingularization is independent of the ambient variety $X$.
Let $\phi_1: U\subset X_1$ and $\phi_2: U\subset X_2$  be two closed embeddings and
let $\widetilde{U}_i\subset \widetilde{X}_i$ be two embedded desingularizations.
Find  embeddings $\psi_i: X_i\to {\bf A}^n$ into affine space ${\bf A}^n$. They define the embeddings 
$\psi_i\phi_i: U\to {\bf A}^n$. By Lemma \ref{le: l}, there are embeddings $\Psi_i: {\bf A}^n\to {\bf A}^{2n}$ such that
$\Psi_1\psi_1\phi_1=\Psi_2\psi_2\phi_2: U\to {\bf A}^{2n}$. Since embedded desingularizations commute with closed embeddings of ambient varieties we see that the $\widetilde{U}_i$ are isomorphic over $U$. 

(2) Let $V\to U$ be an open embedding of affine varieties. Assume first that $V=U_f=U\setminus V(f)$, where $f\in K[U]$ is a regular function on $U$.
Let $U\subset X$ be a closed embedding into an affine variety $X$. Then
$U_f\subset X_F$ is a closed embedding into an affine variety $X_F=X\setminus V(F)$ where $F$ is a regular function on $F$ which restricts to $f$.
Since embedded desingularizations commute with smooth morphisms the open embedding
$X_F\subset X$ defines the open embedding of embedded desingularizations $(\widetilde{X_F},\widetilde{U_f})\subset (\widetilde{X},\widetilde{U})$ and the open embedding of desingularizations $\widetilde{U_f}\subset \widetilde{U}$.

Let $V\subset U$ be any open subset which is an affine variety. Then there are desingularizations $\res_V: \widetilde{V}\to V$ and $\res_U:\widetilde{U}\to U$.
Suppose the natural birational map $\widetilde{V} \to \res_U^{-1}(V)$ is not an isomorphism
over $V$. Then we can find an open subset $U_f\subset V$ such that 
$\res_V^{-1}(U_f) \to \res_U^{-1}(U_f)$ is not an isomorphism over $U_f$. But $U_f=V_f$ and by the previous case $\res_V^{-1}(U_f)\simeq \widetilde{U_f}=\widetilde{V_f}\simeq \res_U^{-1}(V)$.\qed

\bigskip
Let $U_i$ be an open affine cover of $X$. For any two open subsets $U_i$ and $U_j$
set $U_{ij}:=U_i\cap U_j$. For any $U_i$ and $U_{ij}$ we find  canonical resolutions $\widetilde{U_i}$ and $\widetilde{U}_{ij}$ respectively.  We define $\widetilde{X}$ to be a variety obtained by glueing $\widetilde{U_i}$ along $\widetilde{U}_{ij}$. Then $\widetilde{X}$
is a smooth variety and $\widetilde{X}\to X$ defines a canonical desingularization independent of the choice of $U_{ij}$.

\subsection{Commutativity of non-embedded desingularization with smooth morphisms} 
\begin{lemma} Let $\phi_1: U\hookrightarrow {\bf A}^m$ and $\phi_2: V\hookrightarrow {\bf A}^n$  be  closed embeddings of  affine varietes $U$ and $V$.  Let $\phi: U\to V$ be an \'etale morphism at $0\in U$.
Then there exists a variety $X\subset {\bf A}^m$ containing $U$ and smooth at $0$ and a morphism $\Phi: X\to {\bf A}^n$ extending the  morphism $\phi: U\to V$ and which is  \'etale at $0$.
\end{lemma}
\noindent{\bf Proof.} Let $\overline{x}:=x_1,\ldots, x_n$ and $\overline{y}:=y_1,\ldots,y_m$  be coordinates on ${\bf A}^n$ and ${\bf A}^m$ respectively. Let $g_1:=\phi^*(x_1),\ldots,g_n:=\phi^*(x_n)$ be generators of the ring
$K[V]\subset K[U]$. Write $K[V]=K[x_1,\ldots,x_n]/(f_1(\overline{x}),\ldots,f_l(\overline{x}))$. Extending  the set of generators of $K[V]$ to a set of generators of $K[U]$ gives   
$$K[U]=K[x_1,\ldots,x_n,y_1,\ldots,y_m]/(f_1(\overline{x})\ldots, f_l(\overline{x}),h_1(\overline{x},\overline{y}),\ldots,h_r(\overline{x},\overline{y})).$$
Since $\phi$ is \'etale at $0$ the functions $x_1,\ldots,x_n$ generate the maximal ideal of $$\widehat{\cO}_{0,U}=K[[x_1,\ldots,x_n,y_1,\ldots,y_m]]/(f_1(x),\ldots,f_l, h_1,\ldots,h_r)=K[[x_1,\ldots,x_n]]/(f_1,\ldots,f_l)=\widehat{\cO}_{x,V}.$$ Choose a maximal subset $\{h_{i_1},\ldots,h_{i_s}\}\subset \{h_1,\ldots,h_r\}$ for which $x_1,\ldots x_n,h_{i_1},\ldots,h_{i_s}\in \frac{(x_1,\ldots x_n,y_1,\ldots, y_m)}{(x_1,\ldots x_n,y_1,\ldots, y_m)^2}$ are linearly independent. Then $s=m$ and $(x_1,\ldots x_n,h_{i_1},\ldots,h_{i_m})=(x_1,\ldots,x_n,y_1,\ldots,y_m)$ define the set of parameters at $x$. 

The subvariety 
$X=\{(x,y)\mid h_{i_1}=\ldots=h_{i_m}=0\}\subset {\bf A}^m$ is smooth at $0$ and the restriction $p_{|X}: X\to {\bf A}^n$ of the natural projection $p:{\bf A}^{m+n}\to {\bf A}^{n}$  is \'etale at $0$. Consequently, 
$U'\to V$ is \'etale at $0$ where $U':=\Spec(K[x_1,\ldots,x_n,y_1,\ldots,y_m]/(f_1(x),\ldots f_l(x), h_{i_1}(x,y),\ldots,h_{i_m}(x,y))$. Also  the closed embedding $U\to U'$ is \'etale at $0$. Then $U$ is a component of $U'$.\qed

Every smooth morphism $\phi: U\to V$ of relative dimension $r$ locally (for some open $U'\subset U$) factors through an \'etale morphism $\psi: U'\to V\times {\bf A}^r$ followed by the natural projection $p:V\times {\bf A}^r\to V$.
By the above and the lemma for any smooth morphism $\phi: U\to V$ and $x\in U$ we can find a neighborhood $U_x\subset U$ of $x$ and a smooth morphism of embedded varieties $(U_x,X_{U_x})\to (V,X_V)$ where $X_{U_x}$ and $X_V$ are smooth. Commutativity of embedded desingularizations with \'etale morphisms implies commutativity of nonembedded desingularizations.

\end{document}